\documentclass[12pt,a4paper]{article}
\usepackage{amsmath}
\usepackage{amssymb}
\usepackage{amsthm}
\usepackage{amsfonts}
\usepackage{latexsym}

\theoremstyle{plain}
\newtheorem{thm}{Theorem}[section]
\newtheorem{cor}{Corollary}[section]
\newtheorem{lem}{Lemma}[section]
\newtheorem{prop}{Proposition}[section]
\theoremstyle{definition}
\newtheorem{defn}{Definition}[section]
\newtheorem{exmp}{Example}[section]
\newtheorem{prob}{Problem}[section]
\theoremstyle{remark}
\newtheorem{claim}{Claim}[section]
\newtheorem{rem}{Remark}[section]

\title{Equivariant asymptotics for Bohr-Sommerfeld \\
Lagrangian submanifolds}
\author{Marco Debernardi $^{(1)}$ and Roberto Paoletti $^{(2)}$
\footnote{\noindent{\bf Address.} $^{(1)}$: Dipartimento di
Matematica \textit{F. Casorati}, Via Ferrata 1, Universit\`{a} di
Pavia, 27100 Pavia, Italy; $^{(2)}$: Dipartimento di Matematica e
Applicazioni, Universit\`a degli Studi di Milano Bicocca, Via R.
Cozzi 53, Edificio U5, 20126 Milano, Italy; {\bf e-mail}:
marco.debernardi@unipv.it, roberto.paoletti@unimib.it }}
\date{}
\begin{document}
\maketitle
\section{Introduction}

Let $M$ be an $\mathrm{n}$-dimensional complex projective
manifold, and let $\Omega$ be a Hodge form on it. Let $(L,h)$ be
an Hermitian ample line bundle on $M$, such that $-2\pi i\,\Omega$
is the curvature of its unique compatible connection. Let
$L^*\supseteq X\stackrel{\pi}{\rightarrow}M$ be the unit circle
bundle in the dual line bundle, and denote by $\alpha \in \Omega
^1(X)$ the normalized connection 1-form on $X$. Thus $\alpha$ is a
contact form on $X$ satisfying $d\alpha =\pi ^*(\Omega)$.

The compact Legendrian submanifolds of $X$ play an important role
in the theory of geometric quantization. An immersed Lagrangian
submanifold $\iota: \Lambda \rightarrow M$ lifts to an immersed
Legendrian submanifold $\tilde \iota: \Lambda \rightarrow X$ if
and only if there exists a non vanishing covariantly constant
section of $\iota ^*(L)$. Thus, the Legendrian submanifolds of $X$
determine by projection distinguished immersed Lagrangian
submanifolds of $(M,\Omega)$, so-called Bohr-Sommerfeld Lagrangian
submanifolds. Roughly speaking, from a semiclassical point of view
these (rather than points in the phase space $(M,\Omega)$) are the
geometric counterparts to physical states. Consequently, the
quantization of Bohr-Sommerfeld Lagrangian submanifolds has been
an important line of research in symplectic geometry (see for
example \cite{gs-gc}, \cite{bpu}, \cite{gt}, \cite{bw} and
references therein).

In particular, a systematic procedure for quantizing
Bohr-Sommerfeld Lagrangian submanifolds has been developed by
Borthwick, Paul and Uribe in \cite{bpu}. In short, the choice of a
half-form $\lambda$ on $\Lambda$ determines a generalized
half-form on $X$ supported on $\Lambda$, essentially the
delta-function determined by $(\Lambda,\lambda)$; by applying the
Szeg\"{o} kernel to the latter, and taking Fourier components, one
then naturally associates to $(\Lambda,\lambda)$ a sequence of
holomorphic sections of $L^{\otimes k}$, $u_k\in H^0(M,L^{\otimes
k})$, for every $k=0,1,2,\ldots$. The theory of \cite{bpu}
describes how the local geometry of $(\Lambda,\lambda)$ captures
the pointwise asymptotic properties of the sequence $u_k$.

In this article, we shall suppose given in addition the
holomorphic action of a $ \mathrm{g}$-dimensional connected
compact Lie group $G$ on $M$, Hamiltonian with respect to
$\Omega$. We shall assume that $0\in \frak{g}^*$ is a regular
value for the moment map $\Phi :M\rightarrow \frak{g}^*$; here
$\mathfrak{g}$ denotes the Lie algebra of $G$. We shall also
suppose that $L$ is an ample $G$-line bundle and that the
Hermitian metric $h$ on $L$ is $G$-invariant.

We recall that, up to topological obstructions, the existence of a
linearization amounts to the existence of a moment map
(\cite{kost}, \cite{gs-gq}, \S 3, and \cite{ggk}, chapter VI).
More precisely, in the presence of a linearization one recovers a
moment map by pairing the connection form on $X$ with the
infinitesimal action of the Lie algebra $\frak{g}$. Conversely, to
a moment map there is associated an infinitesimal action of
$\frak{g}$ on $L$; more precisely, $\xi \in \mathfrak{g}$ acts on
sections of $L$ by the operator $ \nabla _{\xi _M}+2\pi i\Phi
^\xi$, where $\nabla$ is the covariant derivative associated to
the connection, $\xi _M$ is the vector field on $M$ generated by
$\xi$, and $\Phi ^\xi =:\left <\Phi,\xi\right>:M\rightarrow
\mathbb{R}$. The obstruction to extend the infinitesimal action of
$\mathfrak{g}$ to an action of $G$ is of topological nature, and
the extension certainly exists if $G$ is simply connected.

In this situation, every space of global holomorphic sections
$H^0(M,L^{\otimes k})$ admits a $G$-equivariant unitary
decomposition over the irreducible representations of $G$:
\begin{equation*}
H^0(M,L^{\otimes k})=\bigoplus _\varpi H^0(M,L^{\otimes
k})_\varpi.\end{equation*} Here, $\varpi$ runs over the set of
highest weights, and thus indexes all finite dimensional
irreducible representations $V_\varpi$ of $G$; for every $\varpi$,
the summand $H^0(M,L^{\otimes k})_\varpi$ is $G$-equivariantly
isomorphic to a direct sum of finitely many copies of $V_\varpi$.

In particular, if $u_k\in H^0(M,L^{\otimes k})$ is the sequence
associated to the pair $(\Lambda,\lambda)$, we have for every
$k=0,1,2,\ldots$ a decomposition $u_k=\oplus _\varpi
u_{k,\varpi}$, where $u_{k,\varpi}\in H^0(M,L^{\otimes
k})_\varpi$. We shall investigate the asymptotic properties of the
sequence $u_{k,\varpi}$, for $\varpi$ fixed and $k\rightarrow
+\infty$. Naturally enough, these are governed by the mutual
position of $\Lambda$ and the zero locus of the moment map, $\Phi
^{-1}(0)\subseteq M$.

For example, when $G$ is semi-simple and $\Lambda$ covers a
Lagrangian submanifold $\pi (\Lambda)\subseteq M$, $\Lambda$ is
$G$-invariant if and only if $\pi(\Lambda) \subseteq \Phi
^{-1}(0)$ \cite{gs-gq}. If we choose, as we may after averaging,
$\lambda$ itself to be $G$-invariant, then so will be each $u_k$;
therefore, $u_{k,\varpi}=0$ for every $\varpi \neq 0$ and $k\in
\mathbb{N}$.

We shall assume instead that $\Lambda$ is \textit{transversal} to
$\Phi ^{-1}(0)$; this geometric hypothesis implies a nontrivial
decomposition over the irreducibles of $G$.

Incidentally, we remark that any given compact Legendrian
submanifold $\Lambda\subseteq X$ may be deformed into one
transversal to $\Phi ^{-1}(0)$, by a contactomorphism arbitrarily
close to the identity. To see this, for some integer $r\ge 1$ let
us choose Hamiltonian vector fields $V_1,\ldots, V_r$ on $M$, such
that for every $m$ in an open neighbourhood $U$ of $\pi (\Lambda)$
one has $T_mM=\mathrm{span}\{V_1(m),\ldots,V_r(m)\}$. For every
$i=1,\ldots, r$, let $\psi _i:\mathbb{R}\rightarrow
\mathrm{Diff}(M)$ be the one-parameter group of symplectomorphisms
generated by $V_i$, and consider the smooth map $\Psi :
\Lambda\times \mathbb{R}^r\rightarrow M$ given by $\Psi \big
(x,\overrightarrow{t} \big )=\psi _1(t_1)\circ \cdots \circ\psi
_r(t_r)\circ \pi(x)$ ($x\in \Lambda$,
$\overrightarrow{t}=(t_1,\ldots,t_r)\in \mathbb{R}^r$). By the
assumption on the $V_i$'s, we can find $\delta
>0$ such that, if $B_\delta(0)\subseteq \mathbb{R}^r$ is the open
ball centered at the origin of radius $\delta$, the restriction of
$\Psi$ to $\Lambda \times B_\delta(0)$ is a submersion. The
transversality theorem \cite{gp} then implies that we can find
arbitrarily small $\overrightarrow{t}\in B_\delta(0)$ such that
the map $\psi _1(t_1)\circ \cdots \circ\psi _r(t_r)\circ \pi
:\Lambda\rightarrow M$ is transversal to $\Phi^{-1}(0)$. For every
$i=1,\ldots,r$, there exist vector fields $\widetilde{V_i}$ on $X$
which are $\pi$-related to the $V_i$'s (i.e., the horizontal
component of $ \widetilde{V}_i$ is the horizontal lift of $V_i$ to
$X$, for every $i=1,\ldots,r$), and which generate a one-parameter
group of contactomorphisms
$\widetilde{\psi}_i:\mathbb{R}\rightarrow \mathrm{Diff}(X)$
(\cite{wein}, \S 4, and \cite{gei}, Theorem 2.2). Since
$\widetilde{\psi}_i(t)$ covers $\psi _i(t)$ for every $i$ and
$t\in \mathbb{R}$, $\Lambda _{\overrightarrow{t}}=:
\widetilde{\psi} _1(t_1)\circ \cdots \circ\widetilde{\psi} _r(t_r)
\big (\Lambda\big )$ is a Legendrian submanifold of $X$
transversal to $\Phi ^{-1}(0)$.

Here is an explicit example:
\begin{exmp}
\label{exmp:P1-first} Endow $\mathbb{P}^1$ with the Fubini-Study
metric, so that $L$ is the hyperplane bundle. Then $X$ is the unit
sphere $S^3\subseteq \mathbb{C}^2$, with projection $\pi
:S^3\rightarrow \mathbb{P}^1$ given by the Hopf map. Fix
$e^{ia}\in S^1$ ($-\pi\le a\le \pi$). Let $\iota _a
:S^1\rightarrow \mathbb{C}\oplus \mathbb{C}$ be given by $\iota _a
(e^{i\theta})=(\cos (\theta),e^{ia}\,\sin(\theta))$. Then $\iota
(S^1)\subseteq S^3$ is a Legendrian knot for the standard contact
structure, and may therefore be viewed as a Bohr-Sommerfeld
immersed Lagrangian submanifold of $\mathbb{P}^1$. Let us consider
the Hamiltonian action of $S^1$ on $\mathbb{P}^1$ given by
$t\diamond [z_0:z_1]=:[t\,z_0:t^{-1}\,z_1]$, with moment map $\Phi
([z_0:z_1])=(|z_0|^2-|z_1|^2)/\|z\|^2$ (we use $\diamond$ rather
than $\cdot$ to distinguish the action from the ordinary one given
by scalar multiplication). In affine coordinates, $\iota _a(S^1)$
covers the line through the origin of slope $\tan (a)$, and
$\Phi^{-1}(0)$ is the unit circle centered at the origin. This
example may be generalized in any dimension.
\end{exmp}

One motivation for studying this problem comes from the following
natural question: Let us set
\begin{eqnarray}\label{eqn:Mprime} M'=:\Phi ^{-1}(0)\subseteq
M&\mathrm{and}&M_0=:M'/G.\end{eqnarray} Thus, $M_0$ is the GIT
quotient of $M$ with respect to the action of the complexification
$\tilde G$ of $G$, and $(L,h,\Omega)$ descend to corresponding
orbi-objects $(L_0,h_0,\Omega _0)$ on $M_0$. If we set
\begin{eqnarray}\label{eqn:Xprime}X'=:\pi ^{-1}(M')\subseteq X&\mathrm{and}&X_0=:X'/G,
\end{eqnarray} then $X_0$ is the circle orbi-bundle of the Hermitian line orbi-bundle
$(L_0^*,h_0)$. Let us momentarily suppose to fix ideas that $G$
acts freely on $M'$, so that $(M_0,\Omega_0)$ is a K\"{a}hler
manifold, and $L_0$ a honest ample line bundle on it. Now if
$\Lambda \subseteq X$ is a (half-weighted) Legendrian submanifold
transversal to $M'$, the intersection $\Lambda '=\Lambda \cap X'$
determines by projection an immersed (half-weighted) Legendrian
submanifold $\Lambda _0\rightarrow X_0$, which we may think of as
the \textit{reduction} of $\Lambda$. We thus have corresponding
half-forms $u$ on $X$ and $u_0$ on $X_0$ in the images of the
respective Szeg\"{o} projectors; taking Fourier components we
obtain sequences
$$u_k\in  H^0(M,L^{\otimes
k})\,\,\mathrm{and}\,\,(u_0)_k\in H^0(M_0,L_0^{\otimes k}).$$ On
the other hand, it is well-known that for $k=0,1,2,\ldots$ there
is a natural isomorphism $H^0(M,L^{\otimes k})^G\cong
H^0(M_0,L_0^{\otimes k})$ \cite{gs-gq}. One is then led to ask
whether, under the latter isomorphism, $(u_0)_k=u_{k,0}$, at least
in some asymptotic sense. More pictorially, does the principle
\textit{quantization commutes with reduction} also hold for the
single (transverse, semiclassical) state? To leading order, the
relation between the two sequences is governed by the
\textit{effective potential} of the action, defined as the
function $V_{\mathrm{eff}}$ on $M'$ associating to every $p\in M'$
the volume of its $G$-orbit \cite{burnsg}, and a measure of the
mutual position between $\Lambda$ and the $G$-orbit at a given
point (another appearance of the effective potential in
equivariant asymptotics is described in \cite{pao-sq}). Although
$(u_0)_k$ and $u_{k,0}$ have the same order of growth, the answer
to the question above is negative (see Remark
\ref{rem:negative-answer}).

Following Corollary \ref{cor:action-free-case}, we shall also make
some remarks regarding the case where $\Lambda$ is $G$-invariant.

Some general introductory remarks are in order.

First, while we have followed the general philosophy of Borthwick,
Paul and Uribe, we have based our approach on the parametrix for
the Szeg\"{o} kernel constructed by Boutet de Monvel and
Sj\"{o}strand in \cite{bs}, rather than on the theory of
Fourier-Hermite distributions and symplectic spinors as in
\cite{bpu}. This follows the approach to equivariant asymptotics
already used in \cite{pao-mm}, and is inspired by the study of
algebro-geometric Szeg\"{o} kernels by Zelditch in \cite{z} and
its subsequent developments, as in \cite{bsz}, \cite{stz},
\cite{sz} (in \cite{stz}, in particular, the authors work out
scaling asymptotics for toric eigenfunctions). We shall then deal
with half-densities, rather than half-forms, on the given
Legendrian submanifolds.

We have furthermore made extensive use of the notion of Heisenberg
local coordinates introduced in \cite{sz}, for this makes the
relation between the local geometry of $\Lambda$ and the leading
term in the asymptotic expansions particularly explicit and simple
to express.

Thus, even in the action-free case, our proofs and statements
depart somewhat from the corresponding ones in \cite{bpu}.

Finally, we have focused on the case of complex projective
manifolds. However, given the microlocal description of
almost-complex Szeg\"{o} kernels given in \cite{sz}, our arguments
can be extended to the symplectic almost complex category.

Our statements are best expressed by viewing sections of
$L^{\otimes k}$ as equivariant functions on $X$. Given that
$\alpha$ and $\Omega$ endow $X$ with a $G$-invariant volume form,
functions and half-densities on $X$ may be equivariantly and
unitarily identified.

Briefly, let $\mathcal{H}(X)\subseteq \mathcal{C}^\infty(X)$ be
the Hardy space, and let $\mathcal{H}(X)_k$ be the $k$-th isotype
for the $S^1$-action. Then there is a well-known canonical unitary
isomorphism
$$\mathcal{H}(X)_k\cong H^0(M,L^{\otimes k}),$$
and we shall use the same symbol for the holomorphic sections and
the corresponding equivariant functions. Now, if $\Lambda\subseteq
X$ is a compact Legendrian submanifold, the choice of a smooth
half-density $\lambda$ on it determines a generalized half-density
$\delta _{\Lambda,\lambda}$ on $X$ (\S
\ref{subsectn:half-weighted}). Applying the Szeg\"{o} projector to
$\delta _{\Lambda,\lambda}$, and taking Fourier components, we
obtain as before equivariant functions $u_{k,\varpi}$, for every
integer $k$ and highest weight $\varpi$. Our key result concerns
the asymptotic expansion for an appropriate scaling limit of the
sequence $u_{k,\varpi}$.

More precisely, suppose that $x\in X$ and $w\in T_mM$, where
$m=\pi (x)$. We shall often implicitly identify $T_mM$ with the
horizontal tangent space at $x$ determined by the connection,
$H_x(X/M)\subseteq T_xX$. If $x\in \Lambda$, the tangent space
$T_x\Lambda$ may be viewed as a Lagrangian subspace of $T_{m}M$.
Inspired by the results on the scaling limits of Szeg\"{o} kernels
in \cite{bsz} and \cite{sz}, we shall investigate the asymptotic
behavior of $u_{k,\varpi}(x+w/\sqrt{k})$, for $k\rightarrow
+\infty$ and as $w$ varies in $T_mM$. The point $x+w/\sqrt k$ is
only well-defined up to the choice of a coordinate system near
$x$, and the ambiguity is $O(k^{-1})$; the leading order part of
the asymptotic expansion in Theorem \ref{thm:main2} below is then
indipendent of the choice of local coordinates. For concreteness,
we shall at any rate assume that some system of local Heisenberg
coordinates has been fixed.

Let us introduce some further pieces of notation. Here
orthogonality refers to the standard Euclidean structure on
$\mathbb{C}^\mathrm{n}=\mathbb{R}^{\mathrm{n}}\oplus
\mathbb{R}^{\mathrm{n}}$.

\begin{defn}\label{defn:general-defn}
\begin{description}
  \item[i):] If $(h,g)\in S^1\times G$ and $x\in X$, let
$$d_x(h,g):T_xX\rightarrow T_{(h,g)\cdot x}X$$
be the differential of the action.
\item[ii):] As above, set $\Lambda '=:\Lambda\cap X'$. Let
$H(X/M)_{\Lambda'}$ be the restriction  to $\Lambda'$ of the
horizontal tangent bundle; thus, $H(X/M)_{\Lambda'}$ has fiber at
$x\in \Lambda'$ given by $T_{\pi(x)}M$.
%%%
%%%
\item[iii):] For $m\in M$, let $G_m\subseteq G$ be the stabilizer
subgroup of $m$. If $0\in \mathfrak{g}^*$ is a regular value of
$\Phi$, then $G$ acts locally freely on $M'$; therefore, $G_m$ is
a finite subgroup of $G$ for every $m\in M'$ .
%%%
%%%
\item[iv):] If $x\in X$ and $m\in M$ let us denote by
\begin{eqnarray}
\label{eqn:tg-space-to-orbit} \mathfrak{g}_X(x)=T_x(G\cdot
x)\subseteq T_x(X)&\mathrm{and}&\mathfrak{g}_M(m)=T_m(G\cdot
m)\subseteq T_m(M)
\end{eqnarray}
the tangent spaces to the respective $G$-orbits. Now given that
the $G$-action is horizontal on $X'$, for $x\in X'$ we have a
natural identification
\begin{equation}\label{eqn:natural-identification}\mathfrak{g}_X(x)\cong
\mathfrak{g}_M(\pi (x))\subseteq H_x(X/M).\end{equation} It is
well-known that if $m\in M'$ then $\mathfrak{g}_M(m)$ is a
$\mathrm{g}$-dimensional isotropic subspace of $T_{m}M$
\cite{gs-gq}.
%%%%
%%%%%%%%%%%%%
%%%%%%%%%%%%%%%%%%%%
\item[v):] To leading order, we shall see that only the component
of $w$ in a certain n-dimensional real vector subspace $\tilde
N_\Lambda (x)\subseteq T_{\pi (x)}M$ contributes to
$|u_{k,\varpi}(x+w/\sqrt k)|$. More precisely, recall that
$T_x\Lambda$ may be viewed as a Lagrangian subspace of $T_{\pi
(x)}M$. Now if $\Lambda$ is transversal to $X'$, then
$$T_{x}\Lambda \cap \frak{g}_X(x)=\{0\}$$ for every
$x\in \Lambda'$ (Corollary \ref{cor:xi=0}). Thus, if $x\in
\Lambda'$ there is a direct sum decomposition
\begin{eqnarray}\label{eqn:direct-sum-of-tpm}
T_{\pi (x)}M&=&\big (T_x\Lambda+\frak{g}_M(\pi (x))\big )^\perp
\oplus
_{\mathrm{\perp}}\big (T_x\Lambda+\frak{g}_M(\pi (x))\big )\\
&\cong&\left [\big (T_x\Lambda+\frak{g}_M(\pi (x))\big )^\perp
\oplus _{\mathrm{\perp}}\big(T_x\Lambda'^\perp\cap T_x\Lambda\big
)\right ] \oplus \big [T_x\Lambda'\oplus \frak{g}_M(\pi (x))\big
].\nonumber
\end{eqnarray}
If $x\in \Lambda'$, we shall then set
\begin{eqnarray}\label{eqn:tildeN}
\tilde N_\Lambda (x)&=:&\big (T_x\Lambda+\frak{g}_M(\pi (x))\big
)^\perp \oplus _{\mathrm{\perp}}\big(T_x\Lambda'^\perp\cap
T_x\Lambda\big
);\\
 \tilde T_{\Lambda'}(x)&=:&T_x\Lambda'\oplus
\frak{g}_M(\pi (x)).\label{eqn:tildeT}\end{eqnarray} Thus, $T_{\pi
(x) }M=\tilde N_\Lambda (x)\oplus  \tilde T_{\Lambda'}(x)$.  We
shall denote by $\tilde N_\Lambda$ and $\tilde T_{\Lambda'}$ the
rank-n vector sub-bundles of $H(X/M)_{\Lambda'}$ whose fibres at
$x\in \Lambda'$ are given by, respectively, (\ref{eqn:tildeN}) and
(\ref{eqn:tildeT}).
%%%%%%%%%%%%
%%%%%%%%%%%%%%%
%%%%%%%%%%%%%%%%%%%%%
\item[vi):] Suppose again $x\in \Lambda'$, $m=\pi (x)$. Given
$w\in T_mM$, we shall denote by $w_j$, $j=1,2,3,4$, the components
of $w$ in the following intrinsic and unique algebraic
decomposition:

\begin{equation}\label{eqn:intundec}w=w_a+w_b+w_c+w_d,\end{equation} where
$$w_a\in \big (T_m\Lambda+\frak{g}_M(m)\big )^\perp, \,\,
w_b\in (T_m\Lambda ')^\perp \cap T_m\Lambda,\,\,w_c\in
T_m\Lambda',\, \,w_d\in \frak{g}_M(m).
$$
Thus, $w'=:w_a+w_b$ and $w''=:w_c+w_d$ are the components of $w$
in $\tilde N_\Lambda (x)$ and $\tilde T_{\Lambda'}(x)$,
respectively.
%%%%%%%%%%%%%%%%%%%%%%
%%%%%%%%%%%%%%%%%%%%%%%%%%%%
\item[vii):] Let $\mathrm{dens}_\Lambda$ and $\mathrm{dens}
_\Lambda^{(1/2)}$ be the Riemannian density and half-density on
$\Lambda$, respectively. If $\lambda$ is any smooth half-density
on $\Lambda$, we may write $\lambda =f_\lambda \,\mathrm{dens}
_\Lambda^{(1/2)}$ for a unique smooth function $f_\lambda$ on
$\Lambda$.
%%%%%%%%%%%%%%%%%%%%%%%
%%%%%%%%%%%%%%%%%%%%%%%%%%%%%
\end{description}
\end{defn}

The leading order term of the asymptotic expansion for
$u_{k,\varpi}(x+w/\sqrt k)$ will depend on both the effective
volume of the action at $\pi (x)$, and a function expressing a
pointwise measure of the mutual position between the Legendrian
submanifold $\Lambda$ and the $G$-orbit.

Given any $x\in \Lambda'$, let us choose Heisenberg local
coordinates $(\theta,p,q)$ for $X$ centered at $x$. The horizontal
tangent space $H_x(X/M)$ then gets unitarily identified with
$\mathbb{C}^n$, with complex coordinates $z=p+iq$. Perhaps after
applying a unitary transformation in $z$, we may as well assume
that the Lagrangian subspace $T_x\Lambda\subseteq \mathbb{C}^n$ is
defined by $p=0$. Let us choose an orthonormal basis of
$\mathfrak{g}_M(\pi (x))$ (for the induced metric); the inclusion
$\mathfrak{g}_M(\pi (x))\subseteq T_{\pi (x)}M$ is then described
by a linear map
\begin{equation}\label{eqn:inclusion-of-G-orbit}r
\in \mathbb{R}^\mathrm{g}\mapsto R\,r+i\,TR\,r,\end{equation}
where $R$ is an $\mathrm{n}\times \mathrm{g}$ real matrix, $T$ an
$\mathrm{n}\times \mathrm{n}$ real matrix, and they satisfy
\begin{equation}
\label{eqn:conditions-on-R-and-T} \left \{
\begin{array}{rcc}
\mathrm{rank}(R)&=&\mathrm{g}\\
R^tR+R^tT^tTR&=&I_\mathrm{g}.
\end{array}\right.
\end{equation}
Recalling that $\mathfrak{g}_M(m)\subseteq T_mM$ is an isotropic
subspace when $m\in M'$, one can see that the complex matrix
$R^tR+iR^tTR$ is symmetric, has positive definite real part, and
its determinant is independent of the choices involved. Let $\Xi
_{\Lambda}:\Lambda'\rightarrow \mathbb{C}$ be the smooth function
defined by
\begin{equation}\label{eqn:relative-pos-fctn}
\Xi_{\Lambda}(x)=:\frac{\det(R^tR+iR^tTR)^{-1/2}}{V_{\mathrm{eff}}(\pi
(x))}\,\,\,\,\,\,\,\,(x\in \Lambda').\end{equation} The square
root of the determinant is determined according to the conventions
described in \cite{hor}, \S 3.4.

We are now ready to state our main result. Recall that
$\mathrm{n}=\dim _{\mathbb{C}}(M)$,
$\mathrm{g}=\dim_{\mathbb{R}}(G)$.

\begin{thm}
\label{thm:main2} Suppose that $0\in \mathfrak{g}^*$ is a regular
value of $\Phi$, and that the compact Legendrian submanifold
$\Lambda\subseteq X$ is transversal to $\Phi^{-1}(0)$. Let
$\lambda$ be a smooth half-density on $\Lambda$. Fix a highest
weight $\varpi$ for $G$. For $k=0,1,2,\ldots$, let $u_{k,\varpi}$
be the component of $\delta _{\Lambda,\lambda}$ in
$\mathcal{H}(X)_{k,\varpi}\subseteq \mathcal{H}(X)$. Then:
%%%%%%%%%%%%%%%%%%%%%
%%%%%%%%%%%%%%%%%%%%%%
\begin{description}
%%%%%%%%%%%%%%%%%%%%%%
%%%%%%%%%%%%%%%%%%%%%%
  \item[i):] if $x\not\in (S^1\times G)\cdot \Lambda '$,
  then $u_{k,\varpi}(x)=O(k^{-\infty})$ as $k\rightarrow +\infty$;
  \item[ii):] there exist a positive definite metric $S$ on $\tilde N_\Lambda$
  and a real quadratic form $P$ on $H(X/M)_{\Lambda'}$
  such that the following holds. If $x\in (S^1\times G)\cdot \Lambda '$, let
  $(h_j,g_j)\in S^1\times G$, $1\le j\le r_x$, be the finitely many elements such that
  $(h_j,g_j)\cdot x\in \Lambda$. For every $j$, let
  $$x_j=:(h_j,g_j)\cdot x,\,\,\,w_j=[d_{x}(h_j,g_j)](w)\in H_{x_j}(X/M).$$
  For every $w\in T_{\pi (x)}M$, the following asymptotic expansion holds
  for $k\rightarrow +\infty$:
  %%%%%%%%%%%%%%%%%%%%%%%%%%%%%%%%%%%%%%
  %%%%%%%%%%%%%%%%%%%%%%%%%%%%%%%%%%%%%%
  \begin{eqnarray*}
  u_{k,\varpi}(x+w/\sqrt k)&\sim &k^{(\mathrm{n-g})/2}\,\frac{\dim(V_\varpi)}{|G_{\pi (x)}|}\,
  \frac{1}{\pi ^\mathrm{n}}\,\sqrt{\frac{(2\pi)^{(\mathrm{n}+\mathrm{g})}}{2^{\mathrm{g}}}}\,\sum
  _{j=1}^{r_x}\varrho_0^{(j)}(x,\varpi,k,w)\, f_\lambda(x_j)\\
  &&+\sum _{f\ge
  1}k^{(\mathrm{n-g}-f)/2}\,\sum_{j=1}^{r_x}\varrho_f^{(j)}(x,\varpi,k,w),
  \end{eqnarray*}
  %%%%%%%%%%%%%%%%%%%%%%%%%%%%%%%%%%%%
  %%%%%%%%%%%%%%%%%%%%%%%%%%%%%%%%%%%%
  where $G_{\pi (x)}\subseteq G$ is the stabilizer of $\pi(x)$, and
  for every $j=1,\ldots,r_x$ we have
  $$\varrho_0^{(j)}(x,\varpi,k,w)=:h_{j}^{-k}\,\chi _\varpi (g_j)\,\Xi _{\Lambda}
  (x_j)\,e^{-S_{x_j}(w_j',w_j')-i\,P_{x_j}(w_j,w_j)}.
  $$
Here $\chi _\varpi:G\rightarrow \mathbb{C}$ denotes the character
of the representation $V_\varpi$.

\end{description}

\end{thm}

The real quadratic forms $S$ and $P$ will be described precisely
in the course of the proof (see (\ref{eqn:local-normal-gauge}) and
(\ref{eqn:transv-quadr-form})); they are determined by $R$ and
$T$.

Letting $w=0$, we obtain an asymptotic expansion for
$u_{k,\varpi}(x)$:
\begin{eqnarray*}
  u_{k,\varpi}(x)&\sim &k^{(\mathrm{n-g})/2}\,\frac{\dim(V_\varpi)}{|G_{\pi (x)}|}\,
  \frac{1}{\pi ^\mathrm{n}}\,\sqrt{\frac{(2\pi)^{(\mathrm{n}+\mathrm{g})}}{2^{\mathrm{g}}}}\,\sum
  _{j=1}^{r_x}h_{j}^{-k}\,\chi _\varpi (g_j)\,\Xi _{\Lambda}
  (x_j)\, f_\lambda(x_j)\\
  &&+\mathrm{L.O.T.}.
  \end{eqnarray*}

Let us see what the asymptotic expansion of Theorem
\ref{thm:main2} looks like in the action-free case. In this case
obviously $\Lambda'=\Lambda$, and $\varpi$ may be disregarded.
Here we work in a system of adapted Heisenberg local coordinates
(Definition \ref{defn:adapted-heis-coord}). As a corollary of (the
proof of) Theorem \ref{thm:main2}, we obtain (cfr Theorem 3.12 of
\cite{bpu}):

\begin{cor}
\label{cor:action-free-case} Let $\Lambda\subseteq X$ be any
compact Legendrian submanifold, $\lambda$ a smooth half-density on
$\Lambda$. Let $u_k\in \mathcal{H}(X)_k$ be the components of
$\delta _{\Lambda,\lambda}$ in $\mathcal{H}(X)_k$,
$k=0,1,2,\ldots$. Then:
\begin{description}
  \item [i):] if $x\not\in S^1\cdot \Lambda$, then
  $u_k(x)=O(k^{-\infty})$, as $k\rightarrow +\infty$;
  \item [ii):] if $x\in S^1\cdot \Lambda$, let
  $h_1,\ldots,h_{r_x}\in S^1$ be the finitely many elements
  such that $x_j=:h_j\cdot x\in \Lambda$. Set $m=\pi (x)$, $m_j=\pi (x_j)$.
  Suppose that
  $(\theta,p,q)$ is a system of local Heisenberg coordinates for
  $X$
  adapted to $\Lambda$ at $x$.
  Then for every $w\in T_mM$ the following asymptotic expansion holds as $k\rightarrow
  +\infty$:
\begin{eqnarray*}
  u_{k}(x+w/\sqrt k)&\sim &k^{\mathrm{n}/2}\,
  \left (\frac{2}{\pi }\right )^{\mathrm{n}/2}\,\sum
  _{j=1}^{r_x}\varrho_0^{(j)}(x,k,w)\, f_\lambda(x_j)\\
  &&+\sum _{f\ge
  1}k^{(\mathrm{n}-f)/2}\,\sum_{j=1}^{r_x}\varrho_f^{(j)}(x,k,w),
  \end{eqnarray*}
  where for every $j=1,\ldots,r_x$ we have
  $$\varrho_0^{(j)}(x,k,w)=:h_{j}^{-k}\,e^{-\|w_j^{\perp}\|^2-i\,\Omega_{m_j}(w_j^\perp
  ,
  w_j^{\|})};
  $$
  here $w_j^\perp$ is the component of $w_j=d_xh_j(w)$
  perpendicular to $T_{x_j}\Lambda$, and $w_j^{\|}=w_j-w_j^\perp\in
  T_{x_j}\Lambda$.
\end{description}

\end{cor}

In particular, for $w=0$ we obtain:
\begin{eqnarray*}
  u_{k}(x)&\sim &k^{\mathrm{n}/2}\,\left (\frac 2\pi\right )^{\mathrm{n}/2}\,\sum
  _{j=1}^{r_x}h_{j}^{-k}\, f_\lambda(x_j)+\mathrm{L.O.T.}.
  \end{eqnarray*}

As a consequence of Corollary \ref{cor:action-free-case}, let us
momentarily return to the equivariant setting and take up again
the case of an invariant Legendrian submanifold. More precisely,
suppose that $G$ is semi-simple and acts freely on $\Phi^{-1}(0)$.
Suppose also that $\Lambda \subseteq X$ is a $G$-invariant compact
Legendrian submanifold which maps down diffeomorphically under
$\pi$ onto a Lagrangian submanifold $\pi (\Lambda)\subseteq M$.
Thus, $\pi ( \Lambda)\subseteq \Phi ^{-1}(0)$, and if we choose a
$G$-invariant smooth half-density $\lambda =f_\lambda \cdot
\mathrm{dens}^{(1/2)}_\Lambda$ on $\Lambda$, the corresponding
generalized half-density $\delta _{\Lambda,\lambda}$ is also
$G$-invariant. As we have mentioned, we then have $u_{k,\varpi}=0$
unless $\varpi =0$, and $u_{k,0}=u_k$ for every $k$. Now, $\Lambda
_0=:\Lambda /G$ is a compact Legendrian submanifold of $X_0=X'/G$;
let us define (with a slight abuse of language) $\lambda
_0=:f_\lambda \cdot \mathrm{dens}^{(1/2)}_{\Lambda _0}$. Then it
follows from Corollary \ref{cor:action-free-case} that, up to the
multiplicative factor $(2/\pi) ^{\mathrm{g}/2}\,k^{\mathrm{g}/2}$,
the asymptotic expansion for the corresponding sequence
$\widetilde{u}_{k}$ has the same leading order term as the
asymptotic expansion of $u_k$.

Let us illustrate the Theorem and the Corollary with some
examples.

\begin{exmp}\label{exmp:P1}
Let us consider again the setting of Example \ref{exmp:P1-first}.
Thus, $\Lambda \subseteq S^3\subseteq \mathbb{C}^2$ is the
Legendrian knot given by $\iota (e^{i t})=\big(\cos (t),\sin (
t)\big )$. Let us choose the Riemannian half-density on it, so
that $f_\lambda =1$.

Let us consider first the action-free case. The Szeg\"{o} kernel
at the level $k$ is given by
$$
\Pi _k(x,y)=\frac{(k+1)}{\pi}\,\left <x,y\right
>^k\,\,\,\,\,\,(x,y\in S^3),$$
where $ \left <x,y\right >$ denotes the standard Hermitian product
of $x,y\in \mathbb{C}^2$ \cite{bsz}. Since $\Pi_k$ is self-adjoint
with respect to the $L^2$-Hermitian pairing, we have %%
\begin{equation}\label{eqn:conjugate-conjugate}
\left <\Pi _k(\delta _{\Lambda,\lambda}),f\right >=\left< \delta
_{\Lambda,\lambda},\overline{\Pi _k(\overline{f})}\right
>,\end{equation}
for every $f\in \mathcal{C}^\infty (S^3)$ (we are identifying
half-densities with functions by the Riemannian half-density on
$S^3$). Thus, setting $x=(x_0,x_1)\in S^3\subseteq \mathbb{C}^2$,
we have
\begin{eqnarray}\label{eqn:explicit-on-P1}
u_k(x)&=&\int _{0}^{2\pi}\Pi _k\big (x,(\cos (t),\sin(t))\big )\,
dt=\frac{(k+1)}{\pi}\,\int _0^{2\pi}(x_0\, \cos
(t)+x_1\, \sin (t))^k\,dt \nonumber\\
&=&\frac{(k+1)}{\pi}\,\int _0^{2\pi}e^{\,k\,\log\big (x_0\,\cos
(t)+x_1\, \sin (t)\big )}\,dt=\frac{(k+1)}{\pi}\,\int
_0^{2\pi}e^{i\,k\,S(t,x)}\,dt,
\end{eqnarray}
where $S(t,x)=:-i\,\log(x_0\,\cos (t)+x_1\, \sin (t))$ (any branch
of the logarithm may be used). The latter equalities are
meaningless at those values of $t$ where $ x_0\,\cos (t)+x_1\,
\sin (t)=0$; however, the contribution of a neighbourhood of
radius $\epsilon$ of any of these points is $O(k^{\log
(\epsilon)})$. We shall implicitly introduce cut-off functions
vanishing in a small neighbourhood of those points and ignore them
in the following.

Clearly, $ \Im (S)\ge 0$. We have
$$\frac{\partial S }{\partial t}=-i\,\frac{-x_0\,\sin
(t)+x_1\,\cos (t)}{x_0\,\cos (t)+x_1\,\sin (t)}.$$ Therefore,
$\frac{\partial S }{\partial t}(t_0,x)=0$ if and only if there
exists  $e^{ih}\in S^1$ such that $e^{ih}\cdot x=\iota (t_0)$.
Thus, (\ref{eqn:explicit-on-P1}) is rapidly decreasing in $k$
unless $x\in S^1\cdot \Lambda$. Suppose then $x\in S^1\cdot
\Lambda$, and let $e^{ih_j}\in S^1$, where $h_0,\ldots,h_r\in
[0,2\pi)$, be the elements such that $e^{ih_j}\cdot x\in \Lambda$.
For every $j=1,\ldots,r$ there is a unique $t_j\in [0,2\pi)$ such
that $e^{ih_j}\cdot x=\iota (t_j)$. Hence the $t_j$'s are the only
stationary points of $S$. At every $h_j$, we have
$$
\frac{\partial ^2S}{\partial t^2}=i,
$$
so that the $t_j$ are all non-degenerate critical points. We have
$$
x_0\,\cos(t_j)+x_1\,\sin (t_j)=e^{-ih_j}.$$ Application of the
stationary phase Lemma now yields
\begin{eqnarray}\label{eqn:explicit-on-P1-result}
u_k(x)&\sim&\frac{\sqrt{2\pi}}{\pi}\,\sqrt{k}\cdot \sum
_{j=1}^r\,e^{-ikh_j}+ \mathrm{L.O.T.}
\end{eqnarray}
\end{exmp}

\begin{exmp}\label{exmp:P1-with-action}
Let us re-examine Example \ref{exmp:P1} in the presence of the
action $S^1\times \mathbb{P}^1\rightarrow \mathbb{P}^1$ given by
$t\diamond [z_0:z_1]=:[tz_0:t^{-1}z_1]$ which we considered in
Example \ref{exmp:P1-first}. For any $\varpi \in \mathbb{Z}$, we
now have:
\begin{eqnarray}\label{eqn:explicit-on-P1-equiv}
u_{k,\varpi}(x)&=&\frac{(k+1)}{2\pi^2}\,\int _0^{2\pi}\!\int
_0^{2\pi}e^{i\,k\,S(s,t,x)}\,e^{i\varpi s}dt,
\end{eqnarray}
where $S (s,t,x)=:-i\log\big (x_0\,e^{-is}\,\cos
(t)+x_1\,e^{is}\,\sin (t)\big)$. We have
\begin{eqnarray}\label{eqn:der-wrt-t}
\frac{\partial S}{\partial t}&=&-i\cdot \frac{-x_0\,e^{-is}\,\sin
(t)+x_1\,\cos (t)\,e^{is}}{x_0\,e^{-is}\,\cos
(t)+x_1\,e^{is}\,\sin (t)},\\
\frac{\partial S}{\partial s}&=&\frac{-x_0\,e^{-is}\,\cos
(t)+x_1\,\sin (t)\,e^{is}}{x_0\,e^{-is}\,\cos
(t)+x_1\,e^{is}\,\sin (t)}\label{eqn:der-wrt-s} \end{eqnarray}
Thus, $ d_{s,t}S(s_0,t_0,x)=0$ if and only if $e^{ih}\cdot
(e^{is_0}\diamond x)=\iota (t_0)$ for some $e^{ih}\in S^1$ (by
(\ref{eqn:der-wrt-t})) and $\|x_0\|=\|x_1\|$ (by pairing
(\ref{eqn:der-wrt-t}) with (\ref{eqn:der-wrt-s})). Thus,
$u_{k,\varpi}(x)$ is rapidly decreasing unless $x\in (S^1\times
G)\cdot \Lambda'$ (we have $G=S^1$ here).

Now suppose $x\in (S^1\times G)\cdot \Lambda'$, and let
$(e^{ih_j},e^{is_j})\in S^1\times G$ be the finitely many elements
such that $e^{ih_j}\cdot (e^{is_j}\diamond x)\in \Lambda '$, and
for every $j$ let $t_j\in [0,2\pi)$ be uniquely determined by the
condition that $e^{ih_j}\cdot (e^{is_j}\diamond x)=\iota (t_j)$.
The pairs $(s_j,t_j)$ are the only critical points of
$S(\cdot,\cdot,x)$, and for any $j=1,\ldots,r$ the Hessian matrix
of $S$ at $(s_j,t_j)$ is given by
$$H_{(s_j,t_j)}(S)=i\left(%
\begin{array}{cc}
  1 & \pm i \\
  \pm i & 1 \\
\end{array}%
\right).$$ Applying the stationary phase Lemma, we now obtain:
$$
u_{k,\varpi}(x)\sim \frac{1}{\sqrt 2 \pi}\,\sum _je^{i(\varpi
s_j-kh_j)}+\mathrm{L.O.T.}$$ This agrees with Theorem
\ref{thm:main2}. To check this, remark that $|G_m|=2$ and
$V_{\mathrm{eff}}(m)=\pi$,  for every $m\in \Phi ^{-1}(0)$. The
latter equality follows from the fact that every $G$-orbit in
$S^3$ has length $2\pi$, and if an orbit maps to $\Phi^{-1}(0)$
then it doubly covers its image in $\mathbb{P}^1$.
\end{exmp}

As an application, in \S \ref{sctn:herm-prdcts} we shall study the
following problem:

\begin{prob}
Suppose given two compact Legendrian submanifolds, $\Lambda,
\Sigma\subseteq X$, with specified smooth half-densities $\lambda$
and $\sigma$, respectively. Let $u_{k,\varpi},\,v_{k,\varpi}\in
\mathcal{H}(X)_{k,\varpi}$ be the components of $\delta
_{\Lambda,\lambda}$ and $\delta_{\Sigma,\sigma}$, respectively.
How can we relate the asymptotic behavior of the Hermitian
products $(u_{k,\varpi},v_{k,\varpi})$, as $k\rightarrow +\infty$,
to the geometry of $\Lambda$, $\Sigma$ and $\Phi^{-1}(0)$?
\end{prob}

In the action-free case, and in the setting of Fourier-Hermite
distributions and symplectic spinors, this was carried out in
\cite{bpu}.

Broadly speaking, we shall see that:

\begin{itemize}
    \item \textit{if $\Big ((S^1\times G)\cdot \Lambda \Big)\cap \Sigma \cap
\Big ((\Phi \circ \pi)^{-1}(0)\Big )=\emptyset$, then
$(u_{k,\varpi},v_{k,\varpi})=O(k^{-\infty})$ as $k\rightarrow
+\infty$;}
    \item \textit{if, more generally, the map $S^1\times G\times \Lambda\rightarrow
X$ given by the action is transversal to $\Sigma'=\Sigma\cap \Big
((\Phi \circ \pi)^{-1}(0)\Big )$, then there is an asymptotic
expansion
$$
(u_{k,\varpi},v_{k,\varpi})\sim k^{-\mathrm{g}/2}\rho _0+\sum
_{f\ge 1} k^{-(\mathrm{g}+f)/2}\rho _f,$$ where the leading term
$\rho _0$ is described explicitly, and is a sum of terms
corresponding to each $(h_j,g_j)\in S^1\times G$ such that
$\Big((h_j,g_j)\cdot \Lambda\Big )\cap \Sigma'\neq \emptyset$;}
    \item \textit{a similar asymptotic expansion holds when $S^1\times G\times
\Lambda\rightarrow X$ meets $\Sigma'$ nicely; the order of the
leading term depends on the dimension of the inverse image of
$\Sigma'$ in $S^1\times G\times \Lambda$, and the leading
coefficient is determined by certain integrals on this inverse
image.}
\end{itemize}

\bigskip

The present work covers part of the PhD thesis of the first author
at the University of Pavia.

\textbf{Acknowledgements:} We are very grateful to the referee for
suggesting various improvements in presentation, and to Steve
Zelditch for some interesting remarks.

\section{Preliminaries}

In this section we shall collect a number of preliminary technical
results, and begin a more precise description of the microlocal
background of the quantization scheme outlined in the
introduction. In \S \ref{subsectn:equiv-setting} we shall prove
statement i) of Theorem \ref{thm:main2}.

\subsection{The distribution defined by a half-density on a Legendrian submanifold.}
\label{subsectn:half-weighted}

The given complex structure $J$ on $M$ and the unique compatible
connection form $\alpha$ on $X$ determine a Riemannian metric and
a volume form $\mathrm{vol}_X=:\alpha \wedge (d\alpha)^n$. The
generator $\frac{\partial}{\partial \theta}$ of the $S^1$-action
on $X$ spans the vertical tangent bundle $V(X/M)=:\ker
(d\pi)\subseteq TX$:
\begin{equation}\label{eqn:vert-tg-bundle}
V(X/M)=\mathrm{span} \left\{\frac{\partial}{\partial
\theta}\right\},\end{equation} and $TX=V(X/M)\oplus H(X/M)\cong
V(X/M)\oplus \pi^*(TM)$.

Now let $\Lambda \subseteq X$ be a compact Legendrian submanifold,
endowed with the induced Riemannian metric. Suppose $x\in
\Lambda$. In view of the above isomorphism, any basis $\frak{b}$
of $T_x\Lambda$ can be naturally extended to a basis
$\widetilde{\frak{b}}=(\frac{\partial}{\partial
\theta},\frak{b},J_{\pi(x)}\frak{b})$ of $T_xX$, where
$J_{\pi(x)}$ denotes the complex structure at $\pi (x)\in M$. By
construction, $\widetilde{\frak{b}}$ is orthonormal if so is
$\frak{b}$. Thus the map $\frak{b}\mapsto \widetilde{\frak{b}}$
yields an embedding
$\mathrm{Bs}(\Lambda)_{\mathrm{ort}}\hookrightarrow
\mathrm{Bs}(X)_{\mathrm{ort}}$, with obvious equivariance
properties; here $\mathrm{Bs}(\Lambda)_{\mathrm{ort}}$ is the
principal $O(n)$-bundle of orthonormal frames of $\Lambda$, and
$\mathrm{Bs}(X)_{\mathrm{ort}}$ is the principal $O(2n+1)$-bundle
of orthonormal frames of $X$. Given this, half-densities on $X$
restrict to half-densities on $\Lambda$; conversely,
half-densities on $\Lambda$ extend to half-densities for $X$
defined on $\Lambda$.

On the upshot, any choice of a smooth half-density $\lambda$ on
$\Lambda$ determines a generalized half-density on $X$ supported
on $\Lambda$, $\delta_{\Lambda,\lambda}$, as follows: If $\beta$
is a smooth half-density on $X$, let $\beta _\Lambda$ denote the
induced half-density on $\Lambda$; thus the product $\beta
_\Lambda \otimes \lambda$ is a density on $\Lambda$. Then
$$
\delta _{\Lambda,\lambda}(\beta)\,=:\,\int _\Lambda \beta _\Lambda
\otimes  \lambda.$$ Suppose $\beta =b\,\mathrm{dens} _X^{(1/2)}$,
$\lambda =f_\lambda \,\mathrm{dens} _\Lambda^{(1/2)}$ with $b\in
\mathcal{C}^\infty (X)$, $f_\lambda \in \mathcal{C}^\infty
(\Lambda)$ (here $\mathrm{dens} _X^{(1/2)},\,\mathrm{dens}
_\Lambda^{(1/2)}$ denote the Riemannian half-densities on $X$ and
$\Lambda$, respectively). Then $\beta _\Lambda=\left .b\right
|_\Lambda \,\mathrm{dens} _\Lambda^{(1/2)}$, and
\begin{equation}\label{eqn:explicit-pairing}
\delta _{\Lambda,\lambda}(\beta)\,=:\,\int _\Lambda \left (\left
.b\right |_\Lambda \cdot f_\lambda \right )\, \mathrm{dens}
_\Lambda.\end{equation}

\subsection{Adapted local coordinates and the microlocal structure of $\delta
_{\Lambda,\sigma}$.}\label{subsctn:adapted}
%
%In the following we shall be studying the asymptotics of the Hardy
%space component $\Pi (\delta _{\Lambda,\sigma})$, where $\Pi$ the
%Szeg\"{o} kernel of $X$. To this end, we shall apply the
%stationary phase Lemma to the integral representing $\Pi (\delta
%_{\Lambda,\sigma})$. The resulting computations will be made
%easier by selecting
%
%
%Here we shall work at a point $x_0\in \Lambda$ and select in the
%neighborhood of $x_0$ a system of local equations for $\Lambda$
%adapted to the CR and contact structures of $X$. We shall rely on

Heisenberg coordinates for circle bundles are discussed in
\cite{sz}. The context of \cite{sz} is the symplectic and almost
complex category; we recall that the construction of local
Heisenberg coordinates at a given $x_0\in X$  involves the choice
of \textit{preferred local coordinates} and a preferred frame for
$L$ at $m_0=\pi (x_0)\in M$. Although it isn't strictly necessary,
in the present complex projective setting the preferred local
coordinates and frames involved may as well be assumed
holomorphic.
%
%
%
%This choice will be based on the notion of local Heisenberg
%coordinates at a point $x_0\in X$, introduced in \S 1.2 of
%\cite{sz-rssm}.
%
%

In  short, suppose that $(z_1,\ldots,z_n)$ is a system of
preferred local holomorphic coordinates for $M$ at $m_0$, so that
the Hermitian metric satisfies $\left (g-i\Omega)\right
|_{m_0}=\sum _{j=1}^n dz_j\otimes d\overline z_j$. Let $e_L$ be a
preferred local holomorphic frame for $L$ at $p_0$, with dual
frame $e_L^*$, such that $e_L^*(m_0)=x_0$. The associated system
of local Heisenberg coordinates for $X$ centered at $x_0$, $\rho
:U\subseteq (-\pi,\pi)\times \mathbb{C}^n\rightarrow V\subseteq
X$, is
$$\rho (\theta ,z)=e^{i\theta}\,a(z)^{-1/2}e_L^*(z);$$
here $a=:\left \|e_L^*\right \|^2=\left \|e_L\right \|^{-2}$.

Write $z=(z_1,\ldots,z_n)=p+iq$, where $p,\,q\in \mathbb{R}^n$,
and set $p\,dq=:\sum _jp_j\,dq_j$, $q\,dp=:\sum _jq_j\,dp_j$. By
\cite{sz}, \S 1.2, the connection form $\alpha$ has the local
representation
\begin{equation}\label{eqn:conn-heis}
\alpha =d \theta +p\, dq-q\,dp+\beta (p,q),\end{equation} with
$\beta =O(\|z\|^2)$.

\begin{defn} Suppose that $\Lambda\subseteq X$ is a compact Legendrian submanifold,
and that $x_0\in \Lambda$. A system of Heisenberg local
coordinates $(\theta,p,q)$ centered at $x_0$ is called
\textit{adapted} to $\Lambda$ at $x_0$ if $\Lambda$ is tangent to
the submanifold $\{\theta=0,p=0\}$ at
$x_0$.\label{defn:adapted-heis-coord}
\end{defn}

Any system of Heisenberg local coordinates at $x_0$ may be turned
into one adapted to $\Lambda$ at $x_0$ simply by applying a
suitable unitary transformation in the $z$ coordinates.

Suppose the Heisenberg local coordinates $(\theta,p,q)$ are
adapted to $\Lambda$ at $x_0$. Then $\Lambda$ is locally defined
by $\theta= f(q)$, $p=h(q)$, where $(f,h):V\rightarrow
\mathbb{R}\times \mathbb{C}^n$ vanishes to second order at $x_0$.
Thus, $F(\theta,q)=\theta-f(q)$ and $H(p,q)=p-h(q)$ are local
defining functions for $\Lambda$ on $V$. Actually,

\begin{lem} $f$ vanishes to \textit{third} order at the origin.
\label{lem:f-3rd-order}
\end{lem}

\textit{Proof.} By assumption, the restriction of $\alpha$ to
$\Lambda$ vanishes identically; therefore, $df =-h\,
dq+q\,dh-\beta (h,q)$, which vanishes to second order at $q=0$.

\bigskip

The $q$'s may be naturally viewed as local coordinates on
$\Lambda$. Let $D_\Lambda(q)$ be the local coordinate expression
for the Riemannian half-density $\mathrm{dens}_\Lambda^{(1/2)}$ on
$\Lambda$. In view of (\ref{eqn:explicit-pairing}), we conclude:

\begin{lem}\label{lem:local-fourier-form}
Suppose $x_0\in \Lambda$, and choose adapted Heisenberg local
coordinates at $x_0$, defined in an open neighborhood $V\ni x_0$.
Up to a smoothing contribution, the restriction of $\delta
_{\Lambda,\lambda}$ to $\mathcal{C}^\infty _c(V)$ is a Fourier
integral
\begin{equation}\label{eqn:local-fourier-form}
\frac{1}{(2\pi)^{\mathrm{n}+1}}\,\int \int _{\mathbb{R}\times
\mathbb{R}^n}e^{i\cdot (\tau F+\eta \cdot H)}\, f_\lambda
(q)\,D_\Lambda (q) d\tau \,d\eta .
\end{equation}
\end{lem}

Here $(\tau,\eta)\in \mathbb{R}\times \mathbb{R}^n$, $\eta \cdot H
=\sum _k \eta _kH _k$, where $H(p,q)=p-h(q)$, and $D_\Lambda$ is
the local coordinate expression for the Riemannian density of
$\Lambda$ (the $q$'s restrict to a system of local coordinates on
$\Lambda$). By our choices, $D_\Lambda (0)=1$. The factor
$(2\pi)^{-(\mathrm{n}+1)}$ in front of
(\ref{eqn:local-fourier-form}) comes from the fact that
$\mathrm{codim}(\Lambda,X)=\mathrm{n}+1$.

Let $ \{V_j\}$ be an open cover of $X$ such that whenever $V_j\cap
\Lambda \neq \emptyset$ there exist Heisenberg local coordinates
on $V_j$ adapted to $\Lambda$ at some $x_j\in V_j\cap \Lambda$.
Let $\sum _j\rho_j=1$ be a partition of unity subordinate to the
open cover $ \{V_j\}$. Then $\delta _{\Lambda,\lambda}=\sum
_j\delta _{\Lambda,\lambda\,(j)}$, where each $\delta
_{\Lambda,\lambda\,(j)}=:\rho _j\,\delta _\Lambda$ is either a
smoothing operator or a Fourier integral as
(\ref{eqn:local-fourier-form}).

The construction of a system of Heisenberg local coordinates
adapted to $\Lambda$ at some $x\in \Lambda$ may be varied smoothly
with $x$. More precisely, let $B_{2n+1}(0,\varepsilon)\subseteq
\mathbb{R}^{2n+1}\cong \mathbb{R}\times\mathbb{C}^n $ be the open
ball of radius $\varepsilon$ centered at the origin and having
radius $\varepsilon >0$. Then:

\begin{lem}\label{lem:variational-version}
Fix $y\in \Lambda$. Then there exist
\begin{description}
  \item[i):] open neighborhoods $y\in U\subseteq \Lambda$ and $y\in V\subseteq
  X$, and
  \item[ii):] a smooth map $\kappa:U\times
B_{2n+1}(0,\varepsilon)\rightarrow V$,
\end{description}
such that the following holds: For every $x\in U$, the restricted
map $\kappa_x=\kappa(x,\cdot): B_{2n+1}(0,\varepsilon)\rightarrow
V$ is a Heisenberg local chart adapted to $\Lambda$ at $x$.
\end{lem}

Let $(\theta ^{(x)}, p^{(x)}, q^{(x)})$ be the local coordinates
associated to $\kappa_x$ ($x\in U$). It is then clear that we may
also find smoothly varying local defining functions
$(F_x,H_x):V\rightarrow \mathbb{R}\times \mathbb{C}^n$, of the
form $F_x=\theta ^{(x)}-f _x(q ^{(x)})$, $H_x=p ^{(x)}-h_x(q
^{(x)})$.

\begin{lem}\label{lem:horizontal-lift}
Suppose $x_0\in X$ and fix Heisenberg local coordinates
$(\theta,z)$ centered at $x_0$. Let $m_0=:\pi (x_0)\in M$ and for
some $\delta
>0$ consider a smooth path $\gamma :(-\delta,\delta)\rightarrow M$
satisfying $\gamma (0)=m_0$. Let $\tilde \gamma
:(-\delta,\delta)\rightarrow X$ be the unique horizontal lift of
$\gamma$ to $X$ satisfying $\tilde \gamma (0)=x_0$. Then the local
Heiseinberg coordinates of $\tilde \gamma$ are of the form
$(\theta (t),z(t))$, where $z(t)\in \mathbb{C}^n$ are the
(holomorphic) preferred local coordinates of $\gamma (t)$, and
$\theta (t)\in (-\pi,\pi)$ vanishes to third order at $t=0$.
\end{lem}

\textit{Proof.} If $z(t)$ are the preferred coordinates of $\gamma
(t)$, then clearly the Heisenberg coordinates of $\tilde \gamma
(t)$ have the form $(\theta (t),z(t))$ for some smooth real
function $\theta (t)$ vanishing at the origin. Now let
$z'(0)=p_0+iq_0$ be the tangent vector of $\gamma$ at $t=0$
(expressed in local coordinates). Thus,
$z(t)=(tp_0+P(t))+i(tq_0+Q(t))$, where $P(t)$ and $Q(t)$ vanish to
second order at $t=0$. Since $\tilde \gamma$ is horizontal, the
pull-back $\tilde \gamma ^*(\alpha)\in \Omega ^1(-\delta,\delta)$
vanishes identically. Now in view of (\ref{eqn:conn-heis}) and the
horizontality of $\tilde \gamma$,
\begin{eqnarray*}
0=\tilde \gamma ^*(\alpha)&=&\big [\theta
'(t)+(tp_0+P(t))\,(q_0+Q'(t))-(tq_0+Q(t))\,(p_0+P'(t))\big
]\,dt\\
&&+O(t^2).
\end{eqnarray*}
It follows that $\theta '(t)=O(t^2)$, since the first order terms
cancel out.

\subsection{The equivariant setting.}\label{subsectn:equiv-setting}

Recall that the connection form $\alpha$ and $\pi^*(\Omega)$
naturally endow $X$ with a $G$-invariant volume form. This yields
a unitary and equivariant identification of functions and
half-densities, which with some abuse of language will be implicit
in the following discussion.

Let $L^2(X)$ be the Hilbert space of square-integrable
half-densities on $X$. By the theory of \cite{bs}, the Schwartz
kernel of the Szeg\"{o} projector $\Pi _X:L^2(X)\rightarrow
\mathcal{H}(X)\subseteq L^2(X)$ is microlocally a Fourier integral
\begin{equation}\label{eqn:parametrix}
\Pi _X(x,y)\,=\, \int _0^{+\infty}\,e^{it\psi (x,y)}\,\zeta
(x,y,t)\,dt.\end{equation} The phase $\psi$ is the restriction to
$X\times X$ of a smooth function on $L^*\times L^*$ defined in the
neighbourhood of $(x_0,x_0)$, satisfying $\Im (\psi)\ge 0$. The
Taylor series of $\psi$ along the diagonal $\Delta _{L^*}\subseteq
L^*\times L^*$ is:
$$\psi (x+h,x+k)\sim i\sum _{I,J}
\frac{\partial ^{I+J}\rho}{\partial z^I\partial \overline
z^J}(x)\, h^I\,\overline k^J,$$ where $\rho =1-\left \|\cdot
\right \|^2$ is the defining function for $X\subseteq L^*$. The
amplitude $\zeta (x,y,t)\in S^n(X\times X\times \mathbb{R}^+)$ is
a classical symbol of the form
\begin{equation}\label{eqn:classical-symb}
\zeta (x,y,t)\sim \sum _{k=0}^\infty t^{n-k}\, \zeta
_k(x,y).\end{equation} A complete discussion of the almost
analytic geometry involved, together with a description of the
leading term, is in \cite{bs}, \cite{z}, \cite{sz}. It follows in
particular that the wave front of $\Pi _X$ is the closed isotropic
cone
$$
\Sigma \,=\,\{(x,r\alpha _x,x,-r\alpha _x): x\in X,r>0\}\subseteq
\left(T^*X\setminus \{0\}\right )\times \left(T^*X\setminus
\{0\}\right ).$$

By standard basic results on wave fronts \cite{duist}, \cite{hor}
$\Pi_X$ extends to a continuous operator
$\Pi_X:\mathcal{D}'(X)\rightarrow \mathcal{D}'(X)$. Its image is
the space of those distributions all of whose Fourier coefficients
belong to the Hardy space.

In particular, if $u=\Pi _X(\delta _{\Lambda,\lambda})$ then the
wave front of $u$ satisfies $\mathrm{WF}(u)\subseteq \{(x,r\alpha
_x):x\in \Lambda,r>0\}$; its projection in $X$ is the singular
support of $u$, which thus satisfies: $\mathrm{SS}(u)\subseteq
\mathrm{SS}(\delta _{\Lambda,\lambda})\subseteq \Lambda$. If
$x\not \in S^1\cdot \Lambda$, then $u$ is smooth on an
$S^1$-invariant neighbourhood of $x$; hence
$u_k(x)=O(k^{-\infty})$.

Now the $G$-action on $X$ induces a unitary representation of $G$
on $L^2(X)$, given by $(g\cdot f)(x)= \mu _{g^{-1}}^*f(x
)=f(g^{-1}\cdot x )$. Given a highest weight $\varpi$ for $G$, let
$L^2(X)_\varpi\subseteq L^2(X)$ be the subspace of those elements
contained in a finite direct sum of copies of $V_\varpi$. The
orthogonal projector $P_\varpi:L^2(X)\rightarrow L^2(X)_\varpi$ is
given by
\begin{eqnarray}\label{eqn:Pomega}
P_{\varpi}&=&\dim (V_\varpi)\,\int _G\, \chi _\varpi (g^{-1})\,\mu
_{g^{-1}}^*\,dg,
\end{eqnarray}
where $\chi _\varpi$ is the character of the representation
$\varpi$ \cite{dixmier}.

We need to recall some basic facts concerning the the microlocal
structure of (\ref{eqn:Pomega}). To this end, let us remark that
the action $\mu:G\times X\rightarrow X$ naturally induces a
Hamiltonian action (for the canonical symplectic structure)
$\tilde \mu :G\times \left(T^*X\setminus \{0\}\right ) \rightarrow
T^*X\setminus \{0\}$. Let $\Psi:T^*X\setminus \{0\}\rightarrow
\mathfrak{g}^*$ be the associated moment map. Then $P_\varpi$ is a
Fourier integral operator, associated to the Lagrangian
submanifold
\begin{eqnarray}\label{eqn:lag-and-pomega}
\Lambda _0&=:&\left \{(\nu _1,\nu _2)\,:\,\Psi (\nu _1)=0, \, \nu
_2=\tilde \mu (g,\nu _1)\right \}\nonumber \\
&\subseteq &\left(T^*X\setminus \{0\}\right )\times
\left(T^*X\setminus \{0\}\right )
\end{eqnarray}
\cite{gs-hq}. Thus, $P_\varpi$ obviously extends to a bounded
operator $P_\varpi:\mathcal{D}'(X)\rightarrow \mathcal{D}'(X)$,
and the composition $P_\varpi \circ \Pi_X$ is a Fourier integral
operator with complex phase, whose wave front satisfies:

\begin{eqnarray}\label{eqn:wave-front-composed}
\mathrm{WF}(P_\varpi \circ \Pi_X)&=&\left \{(x,r\alpha
_x,y,-r\alpha _y)\,:\,r>0,\Psi (x,r\alpha _x)=0, \, y=\mu (g,x)\right \}\nonumber \\
&=&\left \{(x,r\alpha _x,y,-r\alpha _y)\,:r>0,\,\Phi (x)=0, \,
y=\mu (g,x)\right \}.\end{eqnarray} In the first equality, we have
made use of the $G$-invariance of $\alpha$, and in the second we
have used the equality $\Psi (x,r\alpha _x)=r\Psi (x,\alpha
_x)=r\Phi (x)$ \cite{gs-gq}.

On the upshot,
\begin{eqnarray}\label{eqn:wave-front-projection}
\mathrm{WF}(P_\varpi \circ \Pi_X(\delta
_{\Lambda,\lambda}))&\subseteq &\left \{(x,r\alpha _x)\,:\,x\in
G\cdot \Lambda, \, \Phi (x)=0,\,r>0\right \}.\end{eqnarray}
Therefore, if $x\not \in (S^1\times G)\cdot \Lambda'$, then
$P_\varpi \circ \Pi_X(u)$ is smooth on an $S^1$-invariant
neighborhood of $x$. Given that $u_{k,\varpi}$ is the $k$-th
Fourier component of $P_\varpi \circ \Pi_X(u)$, we obtain:

\begin{prop} \label{prop:rapid-decay}
If $x\not \in (S^1\times G)\cdot \Lambda'$, then
$u_{k,\varpi}(x)=O(k^{-\infty})$ as $k\rightarrow
+\infty$.\end{prop}

Let us now dwell on the geometry of $\Lambda'$. To this end, let
us first recall the following basic fact from \cite{gs-gq}:

\begin{lem}\label{lem:M'-is-coisotropic}
For every $m\in M'$, $T_mM'$ is the symplectic annihilator
$\mathfrak{g}_M(m)^0$ of the isotropic subspace
$\mathfrak{g}_M(m)$. In particular, $T_mM'$ is a co-isotropic
subspace of $T_mM$. \end{lem}

We deduce:

\begin{cor}\label{cor:xi=0}
Suppose that $\Lambda$ is transversal to $X'$, and $x\in
\Lambda'$. Then $$T_x\Lambda \cap T_x(G\cdot x)=0.$$
\end{cor}

\textit{Proof.} Let $m=:\pi (x)\in M'$. Since both $T_x\Lambda$
and $T_x(G\cdot x)=\frak{g}_X(x)$ are horizontal subspaces of
$T_xX$, it is equivalent to prove that $T_x\Lambda \cap
\frak{g}_M(m)=0$, where in the latter equality $T_x\Lambda $ is
identified with the Lagrangian subspace $d_x\pi\big
(T_x\Lambda\big)\subseteq T_mM$. Passing to symplectic
annihilators, we have
\begin{eqnarray*}
\big (T_x\Lambda \cap \frak{g}_M(m)\big )^0&=&T_x\Lambda ^0+
\frak{g}_M(m)^0\\
&=&T_x\Lambda + T_mM'=T_mM,\end{eqnarray*} by Lemma
\ref{lem:M'-is-coisotropic} and the transversality assumption.

\bigskip

By horizontality, this means:

\begin{cor}\label{cor:tranverse-to-orbits}
Let $\Lambda\subseteq X$ be a Legendrian submanifold transversal
to $X'$. If $x\in \Lambda '$, we have
$$
T_x\Lambda \cap \left ( V_x(X/M)\oplus \mathfrak{g}_X (x)\right
)=\{0\}.$$\end{cor}

Since the action of $S^1\times G$ on $X'$ is (locally) free, we
conclude:

\begin{cor} Suppose that $\Lambda\subseteq X$ is a compact Legendrian submanifold
transversal to $X'$, and that $x\in (S^1\times G)\cdot \Lambda'$.
There are then only finitely many $(h_j,g_j)\in S^1\times G$ such
that $(h_j,g_j)\cdot x\in \Lambda$.\label{cor:finitelymany}
\end{cor}

\bigskip

Now suppose $x\in \Lambda '$, and let us choose Heisenberg local
coordinates $(\theta, z)$ adapted to $\Lambda$ centered at $x$,
defined on some open neighborhood $V\ni x$. Let $F(\theta,q)=
\theta-f( q)$ and $H (p,q)=:p-h$ be defining functions for
$\Lambda \cap V$, as in \S \ref{subsctn:adapted}. By construction,
the locus $\theta=0$ is tangent to the horizontal tangent bundle
at the origin, and $d_0f=0$. Given that the $G$-action on $X$ is
horizontal at any $x\in X'$, it follows that $d_xF(\xi (x))=0$. In
view of Corollary \ref{cor:xi=0}, we obtain:

\begin{cor}
There exist an open neighbourhood $E$ of $0\in \mathfrak{g}$ and
$c>0$ such that $\|H \big(\mu_{\exp _G(\xi)}( x)\big)\|\ge
c\|\xi\|$ for every $\xi \in E$.\label{cor:lower-bound-for-h}
\end{cor}

\subsection{A unitary invariant of pairs of Lagrangian subspaces}
\label{subsctn:lag-invariant}

The invariant introduced in this section will be used in \S
\ref{sctn:herm-prdcts}. Let $(V,\Omega_V,J)$ be a unitary vector
space; that is, $V$ is a $2r$-dimensional real vector space,
$\Omega_V$ a linear symplectic structure on $V$, and $J\in
\mathrm{GL}(V)$ is a complex structure compatible with $\Omega_V$.
Leaving $\Omega _V$ and $J$ understood, let $U(V)$ denote its
unitary group, and let $\mathrm{Gr}_\mathrm{Lag}(V)$ be the its
Lagrangian Grassmanian (the manifold parametrizing Lagrangian
vector subspaces of $V$). Given $L,L'\in
\mathrm{Gr}_\mathrm{Lag}(V)$, let $U(V)_{L,L'}\subseteq U(V)$ be
the subset of unitary transformations mapping $L$ onto $L'$.

\begin{defn}\label{defn:unitary-inv-lag}
For every $L$, we let $\imath _J(L,L)=1$. If $c=\dim(L\cap L')<n$,
suppose that $\psi \in U(V)_{L,L'}$ satisfies $\psi (L\cap
L')\subseteq L\cap L'$. Let $\mathcal{B}$ be an orthogonal real
basis of $L$ whose first $c$ vectors lie in $L\cap L'$. The matrix
of $\psi$ in the basis $\mathcal{B}$, viewed as an orthonormal
complex basis of $V$, has a block diagonal form, whose first
$c\times c$ block is a real orthogonal matrix and whose second
$(r-c)\times (r-c)$ block is a unitary matrix $A+iB\in U(r-c)$,
where $A,B\in M_{r-c}(\mathbb{R})$ and $B$ is non-singular. Then
$\imath _J(L,L')=:|\det (B)|$.\end{defn}

For example, if $V=\mathbb{C}^2$ with its standard unitary
structure, and $L,L'\subseteq \mathbb{C}^2$ are two distinct lines
through the origin, then $\imath _J(L,L')=|\sin (\vartheta)|$,
where $\vartheta$ is the angle between $L$ and $L'$.

For every $c=0,\ldots,n$, let $$ D_c=:\{(L,L')\in
\mathrm{Gr}_\mathrm{Lag}(V)\times \mathrm{Gr}_\mathrm{Lag}(V):
\dim (L\cap L')=c\}.$$ We leave it to the reader to check the
following:

\begin{lem}\label{lem:well-defined-and-symm}
$\imath _J :\mathrm{Gr}_\mathrm{Lag}(V)\times
\mathrm{Gr}_\mathrm{Lag}(V)\rightarrow \mathbb{R}^*$ is
well-defined, $U(r)$-invariant (with respect to the action $R\cdot
(L,L')=(RL,RL')$), and symmetric. It is continuous on $D_c$, for
every $c=0,\ldots,n$.
\end{lem}

\section{Proof of Theorem \ref{thm:main2}.}

As before, let $u=:\Pi_X (\delta _{\Lambda,\lambda})\in
\mathcal{H}(X)$, and denote by $u_{k,\varpi}\in
\mathcal{H}_{k,\varpi}(X)$ its $S^1\times G$-equivariant
components.

Suppose $x\in (S^1\times G)\cdot \Lambda'$. Choose local
Heisenberg coordinates $(\theta,z)=(\theta,p,q)$ centered at $x$,
defined on an open neighbourhood $V\ni x$ ($z=p+iq$ and $p,q\in
\mathbb{R}^n$). We shall denote by $x+w$ the point in $V$ having
Heisenberg local coordinates $(0,w)$ ($w\in \mathbb{C}^n$).

Given $w\in \mathbb{C}^n$, let us consider the asymptotics of
$u_{k,\varpi}(x+w/\sqrt{k})$ for $\varpi$ fixed and $k\rightarrow
+\infty$. We have:
\begin{eqnarray}\label{eqn:scaling-limit-grp-integr}
u_{k,\varpi} (x+w/\sqrt{k} ) & = & \frac{\dim
(V_\varpi)}{(2\pi)^{\mathrm{n}+2}}\int _G\int _{-\pi}^\pi \,u\left
(\mu _{g^{-1}}\circ r_{e^{i\vartheta}}
 (x+w/\sqrt{k})\right )\nonumber\\
 && \chi _\varpi (g^{-1})\,e^{-ik\vartheta}\,
dg\,d\vartheta .
\end{eqnarray}
Here $\mu$ and $r$ denote the $G$- and $S^1$-actions on $X$ (we
shall occasionally also use a dot to denote group action on a
given point).

Let $\{(e^{i\vartheta _j},g_j)\}$, $1\le j\le N_m$, be the
finitely many elements of $S^1\times G$ such that
$x_j=:(e^{i\vartheta _j},g_j)\cdot x\in \Lambda$ (Corollary
\ref{cor:finitelymany}); $N_m$ depends only on $m=\pi (x)\in M$.
Since the action of $G$ on $X'$ is locally free, but not
necessarily free, it may happen that $x_j=x_{j'}$ for $j\neq j'$.

We shall now show that, perhaps after disregarding a rapidly
decaying contribution, the integration over $S^1\times G$ may be
localized near the $(e^{i\vartheta _j},g_j)$'s. Using standard
basic facts from the theory of wave fronts \cite{duist},
\cite{hor}, and recalling that we are identifying functions and
densities by means of the Riemannian volume forms, one can prove
the following:

\begin{lem}\label{lem:sing-support-and-group}
For $y\in X$, define the smooth map $\Upsilon_y:S^1\times G
\rightarrow X$ by $\Upsilon_y(h,g)=\mu _{g^{-1}}\circ r_\vartheta
(y)$ ($h\in S^1$, $g\in G$). Then:
\begin{description}
  \item[i):] $\Upsilon _y$ is an immersion, for every $y\in X'$;
  \item[ii):] the pull-back $\Upsilon_y^*( u)$ is a well-defined
  generalized half-density on $S^1\times G$;
\item[iii):] the singular support of $\Upsilon_y^*( u)$ satisfies
$$\mathrm{SS}\left ( \Upsilon_y^*( u)\right )\subseteq
\{(h,g)\in S^1\times G:\mu _{g^{-1}}\circ r_h (y)\in \Lambda
'\}.$$
\end{description}

\end{lem}

Now suppose $\epsilon >0$ is suitably small; similarly, choose
  a suitably small open neighborhood $E$ of the unit $e\in G$. For
  every $j=1,\ldots,N_m$, let
  \begin{center}
  $D_j=:\{e^{i\vartheta}:|\vartheta -\vartheta_j|<\epsilon\}$,
  and $E_j=:g_j^{-1}\cdot E\subseteq G$.
  \end{center}
  Thus, $T_j=:D_j\times
  E_j\subseteq S^1\times G$ is an open neighborhood of $(e^{i\vartheta
  _j},g_j^{-1})$. Let $T_0\subseteq S^1\times G$ be an open subset
  such that $(e^{i\vartheta
  _j},g_j^{-1})\not\in \overline T_0$, for every $j$, and such that
  $S^1\times G=\bigcup _{j=0}^{N_m}T_j$.
  Let $\sum _{j=0}^{N_m}\gamma _j(h,g)=1$ be a partition of unity subordinate
  to the open cover
   $\mathcal{T}=\{T_j\}_{j=0}^{N_m}$ of $S^1\times G$.
We have, with $dh=(2\pi)^{-1}\,d\vartheta$:
\begin{eqnarray}\label{eqn:scaling-limit-grp-integr-partioned}
u_{k,\varpi} (x+w/\sqrt{k} ) & = & \frac{\dim
(V_\varpi)}{(2\pi)^{\mathrm{n}+1}}\,\sum _{j=0}^{N_m}\int _{T_j}
\, \gamma _j(h,g)\, u\left (\mu _{g^{-1}}\circ r_{h}
 (x+w/\sqrt{k})\right )\nonumber\\
 && \chi _\varpi (g^{-1})\,h^{-k}\,
dg\,dh = \sum _j\,u_{k,\varpi} (x+w/\sqrt{k} )_j,
\end{eqnarray}
where $u_{k,\varpi} (x+w/\sqrt{k} )_j$ is defined to be the $j$-th
summand in (\ref{eqn:scaling-limit-grp-integr-partioned}).

\begin{lem}\label{lem:0-th-term}
$u_{k,\varpi} (x+w/\sqrt{k} )_0=O(k^{-\infty})$ as $k\rightarrow
+\infty$.\end{lem}

\textit{Proof.} Perhaps after restricting the open neighborhood
$V$ of $x$, we may assume that
$$ \mathrm{dist}_X\left (\mu
_{g^{-1}}\circ r_{h}
 (y),\Lambda '\right )>\epsilon_1$$
 for some given sufficiently small $\epsilon _1>0$ and every $(h,g)\in T_0$, $y\in
 V$.
 Therefore, as $y\in V$ varies, the generalized functions
 $\psi _0(h,g)\,\Upsilon _y^*(u)\in \mathcal{D}'(S^1\times G)$
 are smooth and have bounded derivatives.
 Taking Fourier components, we deduce that there exist $C_N>0$,
 $N=1,2,\ldots,$ such that
$|u_{k,\varpi}(y)_0|<C_N\,k^{-N}$ for every $y\in V$. Since
$x+w/\sqrt k\in V$ for $k\gg 0$, the statement follows.

\bigskip

Next we shall focus on the asymptotics of each $u_{k,\varpi}
(x+w/\sqrt{k} )_j$, $1\le j\le N_p$. Recalling
(\ref{eqn:scaling-limit-grp-integr}), we have:
\begin{eqnarray}\label{eqn:scaling-limit-1-equiv}
u_{k,\varpi} (x+w/\sqrt{k} )_j & = & \frac{\dim
(V_\varpi)}{(2\pi)^{\mathrm{n}+1}}\,\int _{T_j} \!\int _X
   \tilde \Pi \left (\mu _{g^{-1}}\circ r_h
(x+w/\sqrt{k}),y\right )
     \, \delta _{\Lambda,\lambda}(y)\nonumber \\
     &&\gamma _j(h,g)\,\chi _\varpi (g^{-1})\,h^{-k}\,dy\,dg\,dh.
\end{eqnarray}
For every $j$, set $V_j=:(e^{i\vartheta _j},g_j)\cdot V$. We may
assume that if $(h,g)\in T_j$ and $y\not \in V_j$, then
$$ \mathrm{dist}_X\left (\mu
_{g^{-1}}\circ r_{h}
 (x+w/\sqrt k),y\right )>\epsilon_3$$
 for some $\epsilon _3>0$ and every $k\gg 0$. Given that the
 Szeg\"{o} kernel is smoothing away from the diagonal,
it follows from (\ref{eqn:scaling-limit-1-equiv}) that
\begin{eqnarray}\label{eqn:scaling-limit-1-equiv}
u_{k,\varpi} (x+w/\sqrt{k} )_j & \sim & \frac{\dim
(V_\varpi)}{(2\pi)^{\mathrm{n}+1}}\,\int _{T_j} \!\int _{V_j}
   \tilde \Pi \left (\mu _{g^{-1}}\circ r_h
(x+w/\sqrt{k}),y\right )
     \, \delta _{\Lambda,\lambda}(y)\nonumber \\
     &&\gamma _j(h,g)\,\varrho _j(y)\,\chi _\varpi
     (g^{-1})\,h^{-k}\,dg\,dh\,dy,
\end{eqnarray}
for an appropriate compactly supported bump function $\varrho _j$
on $V_j$, identically equal to one near $x_j$. Here the symbol
$\sim$ means that the difference between the left and right hand
side is $O(k^{-\infty})$.

For every $j$, the local Heisenberg coordinates on $V$ determine
by translation local Heisenberg coordinates on $V_j$ centered at
$x_j$. We may then compose with an appropriate $A_j\in U(n)$ in
the $z$-variable, so as to obtain a system of local Heisenberg
coordinates $(\theta ^{(j)},z^{(j)})$ adapted to $\Lambda$ at
$x_j$, in the sense of \S \ref{subsctn:adapted}. With these
coordinates understood, $\mu _{g_j^{-1}}\circ r_{e^{i\vartheta_j}}
(x+w/\sqrt{k})=x_j+w_j/\sqrt{k}$, where $w_j=A_j(w)$; furthermore,
\begin{eqnarray}\label{eqn:new-coordinates}
\mu _{(g_j^{-1}g)^{-1}}\circ r_{e^{i(\vartheta+\vartheta _j)
}}(x+w/\sqrt k)&=&\mu _{g^{-1}}\circ
r_{e^{i\vartheta}}(x_j+w_j/\sqrt k) \\
&&(g\in E,|\vartheta|<\epsilon).\nonumber \end{eqnarray}

To simplify, when focussing on one $j$ at a time, we shall write
$(\theta,p,q)$ for $(\theta^{(j)}, q^{(j)},p^{(j)})$. Thus,
$F_j(q)=\theta-f_j(q)$ and $H _j(p,q)=p-h_j$ will denote local
defining functions for $\Lambda$ in $V_j$.

Thus,
\begin{eqnarray}\label{eqn:scaling-limit-1-equiv-local-coord}
u_{k,\varpi} (x+w/\sqrt{k} )_j & \sim & \frac{\dim
(V_\varpi)}{(2\pi)^{\mathrm{n}+2}}\,e^{-k\vartheta_j}\,\int
_{E}\int _{-\epsilon}^\epsilon \!\int _{V_j}
   \tilde \Pi \left (\mu _{g^{-1}}\circ r_{e^{i\vartheta }}
     (x_j+w_j/\sqrt k),y\right )
     \nonumber \\
     &&\gamma _j(e^{i(\vartheta_j+\vartheta)},g_j^{-1}g)\,\chi _\varpi
     (g^{-1}g_j)\,e^{-ik\vartheta}\,\varrho _j(y)\, \delta
     _{\Lambda,\lambda}(y)\nonumber \\
     &&\,dg\,d\vartheta\,dy.
\end{eqnarray}

We may assume that the Szeg\"{o} kernel can be represented on each
$V_j\times V_j$ by a Fourier integral as in
(\ref{eqn:parametrix}), and that $\delta _{\Lambda,\lambda}$ is
represented on each $V_j$ by a Fourier integral as in
(\ref{eqn:local-fourier-form}); we shall thus apply
(\ref{eqn:parametrix}) to (\ref{eqn:local-fourier-form}), write
the $k$-th Fourier component of the result as an oscillatory
integral, and study the asymptotics of the latter by the Lemma of
stationary phase \cite{hor}.

Let us fix an orthonormal basis of $\frak{g}$, and identify the
latter with $\mathbb{R}^{\mathrm{g}}$. We may assume that the
exponential map, $\exp _G:\frak{g}\rightarrow G$, induces a
diffeomorphism $ E'=:\exp _G^{-1}(E)\rightarrow E$. Thus, the
linear coordinates on $ E'$ become local coordinates on $E$.

\begin{lem} Let $(s_{j},S_j):E'\subseteq \frak{g}\rightarrow \mathbb{R}\times
\mathbb{C}^n$ be defined by the condition that $\mu
_{e^{-\xi}}(x_{j})$ has adapted Heisenberg local coordinates
$(s_{j}(\xi),S_j (\xi))$. Then $S_j$ is an embedding and $s_j$
vanishes to third order at $0\in \mathfrak{g}$.\label{lem:Sj-sj}
\end{lem}

\textit{Proof.} The first statement holds because the $G$-action
on $\Phi^{-1}(0)\subseteq M$ is locally free, and the second
follows from Lemma \ref{lem:horizontal-lift} since $G$ acts
horizontally on $(\Phi \circ \pi)^{-1}(0)\subseteq X$.

\bigskip

By construction of Heisenberg local coordinates, we have:

\begin{lem}\label{lem:lowerupperbound}
Suppose $y\in V_j$ has adapted Heisenberg local coordinates
$(\theta,z)$. Let $\mathrm{dist}_M$ be the geodesic distance
function on $M$. Then (perhaps after restricting $V_j$ and $E'$):
\begin{equation*}
\frac{1}{2}\|S_{j}(\xi)-z\|\le \mathrm{dist}_M\left (\mu
_{e^{-\xi}}(\pi (x_j)),\pi (y)\right )\le
2\|S_{j}(\xi)-z\|.\end{equation*}
\end{lem}

Let us write $\psi _{j}$ and $\zeta _{j}$ for the phase and
amplitude in (\ref{eqn:parametrix}) on $V_j\times V_j$, and define
\begin{eqnarray}\label{eqn:Psij}
\Psi _{j}(\tau,\eta,t,g,\vartheta,y)&=:& t\psi _{j} (\mu
_{g^{-1}}\circ
     r_{e^{i\vartheta}}
(x_{j}+w_{j}/\sqrt{k}),y) \nonumber\\
&&+\tau F_{j}(y)+\eta \cdot H_{j}(y)-\vartheta .
\end{eqnarray}
Clearly, $\Im (\Psi _{j})=t\,\Im (\psi _{j})\ge 0$.

Performing the change of variables $t\mapsto kt$, $\eta \mapsto
k\eta$, $\tau \mapsto k\tau$, we obtain:
\begin{eqnarray}\label{eqn:scaling-limit-3-equiv-j}
u_{k,\varpi} (x+w/\sqrt{k}
     )_{j}
&\sim&k^{\mathrm{n}+2}\,\frac{\dim
(V_\varpi)}{(2\pi)^{\mathrm{n}+2}}\,e^{-ik\vartheta_j}\,\int
_{\mathbb{R}}\!\int _{\mathbb{R}^n}\int _0^{+\infty}\!\int
_E\!\int _{-\epsilon}^{\epsilon}
     \!\int _{V_{j}}\!e^{ik\Psi _{j}}\,
      \chi _\varpi
     (g^{-1}\,g_j)\nonumber \\
     &&\cdot \gamma _j(e^{i(\vartheta_j+\vartheta)},g_j^{-1}g)\,
     \zeta _{j}(\mu _{g^{-1}}\circ
     r_\vartheta
(x_{j}+w_{j}/\sqrt{k}),y,k t)\, \nonumber \\
&&\cdot \varrho _{j}(y)\, f_\lambda (q)\,D_\Lambda (q)\,d\tau\,
d\eta \,dt\,dg\,d\vartheta \,dy.
      \end{eqnarray}

\begin{rem}\label{rem:cut-off}
Arguing as in the proof of Theorems 2.3.1 and 2.2.2 of
\cite{duist}, an oscillatory integral like
(\ref{eqn:scaling-limit-3-equiv-j}) can be evaluated
asymptotically by implicitly introducing a cut-off in the norm of
$(t,\tau,\eta)$, vanishing for large values of argument (this
justifies the integration by parts in the proof of Lemma
\ref{lem:rapid-decay}).\end{rem}

Let us now split the integration in $dg\,dy$ as follows. Let
$\mathrm{dist}_G$ denote the Riemannian distance function on $G$,
and for $k=1,2,\ldots$ define open sets $A_{jk},\,B_{jk}\subseteq
E\times V_{j}$ by
\begin{eqnarray}\label{eqn:ajk-bjk}
A_{jk}&=:&\{(g,y)\in E\times V_{j}\,:\gamma (g,y)>
k^{-1/3}\},\nonumber \\
B_{jk}&=:&\{(g,y)\in E\times V_{j}\,:\gamma (g,y)< 2k^{-1/3}\},
\end{eqnarray}
where \begin{eqnarray}\gamma (g,y)& =:
&\mathrm{max}\{\mathrm{dist}_G(g,e),\mathrm{dist}_M\left (\mu
_{g^{-1}}(\pi (x)),\pi (y)\right )\}.\end{eqnarray}

Let $a _{jk}+b_{jk}=1$ be a partition of unity on $E\times V_{j}$
subordinate to the open cover $E\times V_{j}=A_{jk}\cup B_{jk}$.
By construction, $a _{jk}$ and $b _{jk}$ may be chosen
$S^1$-invariant. In local coordinates, we may actually assume that
\begin{equation}\label{eqn:part-fnctn-jk}
a _{jk}(\xi,z)=a_{j}(\sqrt[3]{k}\,\xi,\sqrt[3]{k}\,z),\,\,b
_{jk}(\xi,z)=b_{j}(\sqrt[3]{k}\,\xi,\sqrt[3]{k}\,z),\end{equation}
for fixed functions $a_{j}$, $b_{j}$. Then
$$u_{k,\varpi}(x+w/\sqrt{k})_{j}\sim u_{k,\varpi} (x+w/\sqrt{k})_{ja}+u_k
(x+w/\sqrt{k})_{jb},$$ where $u_{k,\varpi} (x+w/\sqrt{k})_{ja}$ is
obtained by multiplying the integrand in
(\ref{eqn:scaling-limit-3-equiv-j}) by $a_{jk}(y)$ - and
integration is thus over $A_{jk}$ - and similarly for
$u_{k,\varpi} (x+w/\sqrt{k})_{jb}$.

\begin{lem} \label{lem:rapid-decay}
$u_{k,\varpi} (x+w/\sqrt{k})_{ja}=O(k^{-N})$, $N=1,2,\ldots$.
\end{lem}

\textit{Proof.} In view of Lemma \ref{lem:lowerupperbound}, for
every $(e^{i\vartheta},\exp_G(\xi))\in S^1\times E$, we have
\begin{eqnarray}\label{eqn:distance-and-bounds}
\mathrm{dist}_X(\mu _{e^{-\xi}}\circ r_\vartheta
(x_j+w_j/\sqrt{k}),y)&\ge &\mathrm{dist}_M\left (\mu
_{e^{-\xi}}\left
(\pi(x_j+w_j/\sqrt{k})\right ),\pi (y)\right )\nonumber \\
&\ge & \frac{1}{2}\, \|z-S_{j}(\xi)\|+O(k^{-1/2}).\end{eqnarray}

Denote by $\mathrm{dist}_X$ the geodesic distance function on $X$.
Fix an open neighborhood $R\ni x_j$ with compact closure,
contained in the chosen chart adapted to $\Lambda$. Let
$$C=:\max _{x'\in \overline R}\{\|d_{x'}H _j\|\}.$$
If $x',x''\in R$ have local coordinates $(\theta',z')$,
$(\theta'',z'')$ then
\begin{equation}\label{eqn:differential-of-h-and-bound}
\|H_j(x')-H_j(x'')\|\le C\, \|z'-z''\|.\end{equation} Here
$\|\cdot\|$ is the standard norm in $\mathbb{C}^n$.

Choose $c\in (0,C)$ satisfying the conclusions of Corollary
\ref{cor:lower-bound-for-h}.

Now let us set
\begin{equation}\label{eqn:first-sub-bound}A_{jk}^{(1)}=:\left \{
(g,y)\in A_{jk}:\mathrm{dist}_M\left (\mu _{g^{-1}}(\pi (x)),\pi
(y)\right )> \frac{c}{10\,C}\, \gamma (g,y) \right
\},\end{equation}
\begin{equation}\label{eqn:second-sub-bound}A_{jk}^{(2)}=:\left \{
(g,y)\in A_{jk}:\mathrm{dist}_M\left (\mu _{g^{-1}}(\pi (x)),\pi
(y)\right )< \frac{c}{5\,C}\, \gamma (g,y) \right
\}.\end{equation} Let $\tau _1+\tau_2=1$ be a partition of unity
of $A_{jk}$ subordinate to the open cover $A_{jk}^{(1)}\cup
A_{jk}^{(2)}= A_{jk}$. We may assume that $\tau_1$ and $\tau_2$
are fixed functions of $(g,y)$, independent of $k$.

Suppose first that $(g,y)\in A_{jk}^{(1)}$. By Lemma
\ref{lem:lowerupperbound} and the hypothesis $\gamma
(g,y)>k^{-1/3}$, this implies
$\|z-S_{j}(\xi)\|>c\,k^{-1/3}/(20\,C)$. Given this,
(\ref{eqn:distance-and-bounds}) implies
\begin{eqnarray}\label{eqn:distance-and-bounds-bis}
\mathrm{dist}_X(\mu _{e^{-\xi}}\circ r_\vartheta
(x_j+w_j/\sqrt{k}),y)&\ge &\frac 13 \|z-S_{j}(\xi)\|\end{eqnarray}
for $k\gg 0$. In view of Corollary 1.3 of \cite{bs}, we deduce
with $g=\exp_G(\xi)$:
\begin{eqnarray}
\left |d_t\Psi _j(\mu _{g^{-1}}\circ  r_\vartheta
(x_j+w_j/\sqrt{k}),y)\right |&= &\left
|\psi_j(\mu _{g^{-1}}\circ r_\vartheta (x_j+w_j/\sqrt{k}),y)\right |\nonumber \\
&\ge &\left
|\Im \psi _j(\mu _{g^{-1}}\circ r_\vartheta (x_j+w_j/\sqrt{k}),y)\right |\nonumber \\
 &\ge &C_1\, \|z-S_{j}(\xi)\|^2,\end{eqnarray} with $C_1>0$ an appropriate
 constant.

The differential operator $L^{(1)}=:\psi ^{-1}\, \frac{\partial
}{\partial t}$ is thus well-defined and smooth on $A_{jk}^{(1)}$,
is positively homogeneous of degree $-1$ in $t$, and satisfies
$$L^{(1)}\Psi _{j}=1,\,\,\,\,\,\,\|L^{(1)}(1,y,g)\|\le
C_2/\|z-S_{j}(\xi)\|^2.$$

If on the other hand $(g,y)\in A_{jk}^{(2)}$ and $k\gg 0$, then
necessarily $\mathrm{dist}_G(g,e)=\gamma (g,y)>k^{-1/3}$. If
$g=\exp _G(\xi)$, we also deduce
\begin{eqnarray*} \|S_{j}(\xi) -z\|&\le&
2\,\mathrm{dist}_M\left (\mu _{g^{-1}}(\pi (x)),\pi
(y)\right )< \frac{2c}{5\,C}\gamma (\exp _G(\xi),y)\\
&=& \frac{2c}{5\,C}\,\mathrm{dist}_G(g,e) \le
\frac{4c}{5C}\|\xi\|.\end{eqnarray*} Thus, if $\xi \in \frak{g}$,
$g=\exp _G(\xi)\in R$, and $c>0$ is as in Corollary
\ref{cor:lower-bound-for-h}, then
\begin{eqnarray}\label{eqn:distance-bounds-on-f}
\|d_\eta \Psi _{j}\|=\|H _j(y)\|&=&\|H _j(y)-H _j(\mu _g(x))+H _j(\mu _g(x))\| \nonumber \\
&\ge &\|H _j(\mu _g(x))\|-\|H _j(\mu _g(x))-H _j(y)\|\nonumber \\
&\ge &c\|\xi\|-C\|S_j(\xi) -z\| \nonumber \\
&\ge &\left (c-\frac 45c\right )\, \|\xi\|\ge C_3\,\|\xi\|\ge
C_3\,\|\xi\|^2,
\end{eqnarray}
since we may assume $\|\xi\|<1/2$ if $g=\exp _G(\xi)\in R$.
Arguing as above, one can then produce a linear first order
differential operator $L^{(2)}$ on $A_{jk}^{(2)}$, positively
homogeneous of degree $-1$ in $\eta$ and with no zero order term,
sushc that $L^{(2)}\Psi_{j} =1$, whence
$L^{(2)}(e^{ik\Psi_{j}})=ike^{ik\Psi_{j}}$, and $\|L^{(2)}\|\le
C_4/\|\xi\|^2$.

Then $L=:\tau _1L^{(1)}+\tau _2L^{(2)}$ is a first order linear
partial differential operator on $A_{j}$, positively homogeneous
of degree $-1$ in $(t,\eta)$, having no zero order term, and
satisfying $L\Psi_{j} =1$. Hence
$L(e^{ik\Psi_{j}})=ike^{ik\Psi_{j}}$, and $\|L\|\le
C'/\|(z,\xi)\|^2$ for some $C'>0$. Let $L^T$ be the transpose
operator (norms and transposes are in the given local
coordinates); then for every $s=1,2,\ldots$ there exists a
constant $C_s>0$ such that
\begin{equation}\label{eqn:bount-on-Pt}
\|(L^T)^s\|\le C_s/\|(z,\xi)\|^{2s}.
\end{equation}
In $L^T$ and its powers only $t$- and $\eta$- derivatives occur,
and the coefficients are functions of $(g,y)$.

Now let us set $S=\mathbb{R}\times \mathbb{R}^n\times
(0,+\infty)\times (-2\epsilon,2\epsilon)$, and
$dX=:d\tau\,d\eta\,dt\,d\vartheta$. Let us write $e^{ik\Psi
_{j}}\,F_j$ for the integrand in the expression for $u_{k,\varpi}
(x+w/\sqrt{k})_{ja}$; the latter is obtained from
(\ref{eqn:scaling-limit-3-equiv-j}) by inserting the additional
factor $a_{jk}$. Then (Remark \ref{rem:cut-off})
\begin{equation}\label{eqn:scaling-limit-3-equiv-j-a}
u_{k,\varpi} (x+w/\sqrt{k}
     )_{ja}
\sim k^{\mathrm{n}+2-s}\,\frac{\dim
(V_\varpi)}{(2\pi)^{\mathrm{n}+2} i^s}\,e^{-ik\vartheta_j}\,\int
_{S}
     \!\int _{A_{jk}}\!e^{ik\Psi _{j}}\,(L^T)^s(F_j)\,dX\,dg\,dy
      .
      \end{equation}
Introducing radial coordinates in the $(z,g)$-variables and
invoking (\ref{eqn:bount-on-Pt}),
\begin{eqnarray}\label{eqn:estimate-on-a_k}
\left |u_{k,\varpi} (x+w/\sqrt{k}
     )_{ja}\right |&\le &D_s\,k^{\mathrm{n}+2-s}
     \!\int _{k^{-1/3}}^\infty\!\,r^{2\mathrm{n}+\mathrm{g}-1-2s}\, dr
      \nonumber \\
&=&D'_s\, k^{(\mathrm{n}-\mathrm{g}-s)/3},\nonumber
      \end{eqnarray}
where $D_s,D'_s$ are appropriate positive constants. This
completes the proof of Lemma \ref{lem:rapid-decay}.

\bigskip

We shall now determine the asymptotics of
$u_{k,\varpi}(x+w/\sqrt{k})_{jb}$. To this end, let us perform the
following change of integration variables:
\begin{equation}
\theta'= \theta,\,\,p'= p-h_j(q),\,\,q'= q.
\label{eqn:coordinate-change}
\end{equation} At the origin, the Jacobian of this transformation
is the identity. In the new coordinates, the phase function
(\ref{eqn:Psij}) becomes
\begin{eqnarray}\label{eqn:Psij-new}
\Psi _{j}(\tau,\eta,t,g,\vartheta,\theta',p',q')&=:& t\psi _{j}
\left (\mu _{g^{-1}}\circ
     r_\vartheta
\left (x_{j}+\frac{w_{j}}{\sqrt{k}}+R\left
(\frac{w_{j}}{\sqrt{k}}\right)
\right),\big (\theta',p'+h_j(q'),q'\big )\right ) \nonumber\\
&&+\tau \big (\theta'-f_j(q')\big )+\eta \cdot p'-\vartheta ,
\end{eqnarray}
where $R:\mathbb{C}^\mathrm{n}\rightarrow \mathbb{C}^\mathrm{n}$
vanishes to second order at the origin. Thus, $R_k=:R\left
(\frac{w_{j}}{\sqrt{k}}\right)=O(k^{-1})$ as $k\rightarrow
+\infty$. In the following, we shall work in the new coordinates
and omit the primes for notational simplicity.

We shall next rescale our coordinates by a factor $k^{-1/2}$, as
follows. First, let us rescale the local coordinates on $G$ in the
neighbourhood $E\ni e$, by writing $\xi =\nu /\sqrt k$; thus on
$E$ we have $g=\exp _G( \xi)=\exp _G(\nu /\sqrt k)$. Let us also
rescale in the same manner the new coordinates
(\ref{eqn:coordinate-change}) centered at $x_{j}$ in the
horizontal direction; more precisely, let us write $(\theta,p,q)=
(\theta,r/\sqrt k,s/\sqrt k)$. Here $\nu\in
\mathbb{R}^{\mathrm{g}}$ (given our choice of an orthonormal basis
of $\frak{g}$) and $r,\, s\in \mathbb{R}^\mathrm{n}$. In
Heisenberg coordinates, $(\theta,r/\sqrt k,s/\sqrt k)$ corresponds
to $(\theta, (r+is)/\sqrt k +h_j(s)/k)$.

By the definition (\ref{eqn:ajk-bjk}) of $B_{jk}$, integration in
$d\nu \,dr\,ds$ takes place over a ball of radius $O(k^{1/6})$ in
$\mathbb{R}^{\mathrm{g}}\times \mathbb{C}^\mathrm{n}$.

Let us express (\ref{eqn:Psij-new}) in rescaled coordinates, and
define
\begin{equation} \label{eqn:Psij-rescaled}
\Psi _{jk}(\tau,\eta,t,\nu,\vartheta,\theta,r,s)\,=:\,\Psi
_{j}\left (\tau,\eta,t,\exp _G\left (\frac{\nu}{\sqrt
k}\right),\vartheta,\theta,\frac{r}{\sqrt k},\frac{s}{\sqrt
k}\right).\end{equation} By Lemma \ref{lem:f-3rd-order}, $f_j$
vanishes to third order at the origin. Thus $f_j(s/\sqrt
k)=k^{-3/2}\,f_{jk}(s)$ for a smooth function $f_{jk}$ vanishing
to third order at the origin. We obtain
\begin{eqnarray} \label{eqn:Psij-rescaled-precise}
\Psi _{jk}&=&t\psi _{j} \left (\mu _{e^{-\nu/\sqrt k}}\circ
     r_\vartheta
(x_{j}+w_{j}/\sqrt{k}+R_k),(\theta,k^{-1/2}\,r+k^{-1}\,h_j(s),k^{-1/2}s
)\right ) \nonumber\\
&&+\tau \,\theta+k^{-1/2}\,\eta \cdot r-\vartheta
+k^{-3/2}\,f_{jk}(s).\end{eqnarray} As usual, we may identify
$w_j$ with a tangent vector in $T_{m_j}M$, where $m_j=:\pi (x_j)$;
clearly, $m_j=g_j\cdot m$. Let $\nu _M$ be the vector field on $M$
generated by $\nu \in \frak{g}$.

\begin{lem}\label{lem:adapted}
The adapted Heisenberg coordinates of $\mu _{e^{-\nu/\sqrt
k}}\circ
     r_\vartheta
(x_{j}+w_{j}/\sqrt{k}+R_k)$ are $$\left (\vartheta -\frac{2}{
k}\,\Omega _{m_j}(\nu _M(m_j),w_j) +Q\left
(\frac{w_{j}}{\sqrt{k}},\frac{\nu}{\sqrt{k}}\right),\frac{1}{\sqrt
k}\,\big (w_j -\nu _M(m_j)\big )+T\left
(\frac{w_{j}}{\sqrt{k}},\frac{\nu}{\sqrt{k}}\right)\right),$$
where $Q,T:\mathbb{C}^\mathrm{n}\times
\mathbb{R}^\mathrm{g}\rightarrow \mathbb{C}^\mathrm{n}$ vanish at
the origin to third and second order, respectively.
\end{lem}

Thus, $Q_k=:Q\left
(\frac{w_{j}}{\sqrt{k}},\frac{\nu}{\sqrt{k}}\right)=O(k^{-3/2})$
and $T_k=T\left
(\frac{w_{j}}{\sqrt{k}},\frac{\nu}{\sqrt{k}}\right)=O(k^{-1})$ as
$k\rightarrow +\infty$ for fixed $\nu$ as $k\rightarrow +\infty$.

\textit{Proof.} Clearly, the preferred holomorphic coordinates of
$\mu _{e^{-\nu/\sqrt k}} (m_j+w_{j}/\sqrt{k}+R_k)$ are
$k^{-1/2}\big (w_j-\nu _M(m_j) \big)+T\left
(\frac{w_{j}}{\sqrt{k}},\frac{\nu}{\sqrt{k}}\right)$, for some
$\mathbb{C}^n$-valued function $T$ vanishing to second order at
the origin. By construction, the Heisenberg coordinates of $\mu
_{e^{-\nu/\sqrt k}}\circ
     r_\vartheta
(x_{j}+w_{j}/\sqrt{k}+R_k)$ then have the form $$\left (\theta
(1/\sqrt k),\frac{1}{\sqrt k}\big (w_j-\nu _M(m_j) \big)+T\left
(\frac{w_{j}}{\sqrt{k}},\frac{\nu}{\sqrt{k}}\right)\right),$$ for
some smooth function $\theta :(-\delta,\delta)\rightarrow
\mathbb{R}$.

To determine the latter, let us momentarily set $\vartheta =0$ and
consider the path, defined for sufficiently small $s,t\in
\mathbb{R}$,
$$
\gamma _s:t\in (-\delta,\delta)\mapsto \mu _{e^{-t\,\nu}}
(x_{j}+s\,w_{j}+R(s\,w_j)),$$ where $R$ is as in
(\ref{eqn:Psij-new}); thus, $R(s\,w_j)$ is a smooth function
$(-\delta,\delta)\rightarrow \mathbb{C}^\mathrm{n}$ vanishing to
second order at $s=0$. Let us write $w_j=p_w+iq_w$, $\nu
_M(p_j)=p_\nu+iq_\nu$. The preferred coordinates of $\pi(\gamma
_s(t))$ are given by $(sp_w-tp_\nu)+i(sq_w-tq_\nu)+Q(sw,t\nu)$,
where $Q:\mathbb{C}^\mathrm{n}\times
\mathbb{R}^\mathrm{g}\rightarrow \mathbb{C}^\mathrm{n}$ vanishes
to second order at the origin. Thus, the Heisenberg coordinates of
$\gamma _s(t)$ have the form
$$
\big (\theta (s,t),(sp_w-tp_\nu)+i(sq_w-tq_\nu)+Q(sw,t\nu)\big
),$$ for some real-valued smooth function $\theta (s,t)$.

\begin{claim}
$\theta (s,t)=(s\, t)\cdot d _0+\theta _1(sw,t\nu)$, where $d _0$
depends only on $w$ and $\nu$, while $\theta _1$ vanishes to third
order at $s=t=0$.\label{claim:st}\end{claim}

\textit{Proof.} By construction, $\theta (s,0)$ vanishes
identically. Therefore, $\theta (s,t)= t\,\theta _1(s,t)$, for a
smooth function $\theta_1$. Given that $G$ acts horizontally on
$X'$, Lemma \ref{lem:Sj-sj} implies that $\theta (0,t)$ vanishes
to third order at $t=0$. Thus, $\theta
_1(s,t)=at^2+t^3\,b(t)+s\,d(s,t)$. The claim follows by writing
$d(s,t)=d_0+d_1(s,t)$, where $d_1(0,0)=0$.

\bigskip

Next we shall determine $d_0$ by use of (\ref{eqn:conn-heis}). The
pull-back $\gamma _s^*(\alpha)\in \Omega ^1(-\delta,\delta)$ is
given by
\begin{eqnarray*}
\gamma _s^*(\alpha)(t)&=&\left \{sd_0+\left [
(sq_w-tq_\nu)\,p_\nu-(sp_w-tp_\nu)\, q_\nu \right ]\right \}\,dt+G_1(s,t)\,dt\\
&=&s\,\left [d_0+(p_\nu\, q_w-q_\nu\,p_w )\right
]\,dt+G_1(s,t)\,dt\\
&=&s\,\left [d_0+\Omega _{m_j}(\nu _M(m_j),w_j)\right
]\,dt+G_1(s,t)\,dt,
\end{eqnarray*}
where $G_1(s,t)$ vanishes to second order at $s=t=0$. On the other
hand, if $\nu _X$ is the vector field on $X$ generated by $\nu \in
\frak{g}$, then $$\alpha (\nu _X)=\Phi ^\nu=:<\Phi,\nu>$$
(equation (5.1) of \cite{gs-gq}). Since $\Phi ^{-1}(0)$ is
$G$-invariant, $\left .\frac{\partial \Phi\circ \gamma}{\partial
t}\right |_{(0,t)}=0$ for every $t$. Thus, $$\Phi (\gamma
_s(t))=s\,d_{m_j}\Phi (w_j)+G_2(s,t),$$ where $G_2$ vanishes to
second order at the origin. Hence,
\begin{eqnarray*}
\gamma _s^*(\alpha)(t)&=&-\left <\Phi (\gamma _s(t)),\nu \right
>\,dt\,=-\,\left [s\,\left <d_{m_j}\Phi (w_j),\nu\right >+G_3(s,t)\right ]\,dt \\
&=&-\left [s\,d_{m_j}\Phi^\nu (w_j)+G_3(s,t)\right ]\,dt=\,\left [
-s\Omega _{m_j}(\nu _M(m_j),w_j)+G_3(s,t)\right ]\,dt,
\end{eqnarray*}
where $G_3$ vanishes to second order, and $\Phi ^\nu=:\left
<\Phi,\nu\right
>$ is the Hamiltonian function associated to $\nu$. Thus
$G_1=G_3$, and $d_0=-2\,\Omega _{m_j}(\nu _M(m_j),w_j)$. Lemma
\ref{lem:adapted} now follows in the case $\vartheta =0$ by
letting $s=t=k^{-1/2}$; in general we need only notice that
$r_\vartheta$ corresponds in local Heisenberg coordinates to a
translation by $\vartheta$.

\bigskip
Following the notation of \cite{bsz}, let us set
\begin{eqnarray}\label{eqn:psi-2-bsz}
K _2(u,v)&=:&u\cdot \overline v-\frac
12(\|u\|^2+\|v\|^2)\nonumber \\
&=&i\,\Im (u\cdot \overline v)-\frac 12\|u-v\|^2
\,\,\,\,\,\,(u,v\in \mathbb{C}^\mathrm{n}).\end{eqnarray} In view
of the asymptotic expansion for the phase discussed in the proof
of the scaling limit of the Szeg\"{o} kernel in \cite{bsz} and
\cite{sz}, Lemma \ref{lem:adapted} implies that $\Psi _{jk}$ in
(\ref{eqn:Psij-rescaled-precise}) has an asymptotic $k$-expansion
of the form
\begin{eqnarray}\label{eqn:psi-jk-asympt}
\Psi _{jk}&\sim&i\,t\,\left [1-e^{i\,(\vartheta -\frac{2}{k}\Omega
_{m_j}(\nu _M(m_j),w_j)-\theta)}\right ]-\vartheta
+\tau\,\theta+\frac{1}{\sqrt k}\,\eta \cdot r  \\
&&-\frac{it}{k}\,e^{i(\vartheta-\theta)}\,K _2\left (w_j-\nu
_M(m_j) ,\,r+i\,s\right
)+O(k^{-3/2})\nonumber \\
&=&i\,t\,\left [1-e^{i\,(\vartheta -\theta)}\right ]-\vartheta
+\tau\,\theta+\frac{1}{\sqrt k}\,\eta \cdot r \nonumber \\
&&-\frac{t}{k}\,e^{i(\vartheta-\theta)}\,\left [2\,\Omega
_{m_j}(\nu _M(m_j),w_j)+iK _2\left (w_j-\nu_M(m_j),\,r+i\,s\right
)\right
]\nonumber \\
&&+t\,e^{i(\vartheta-\theta)}\,P\left ( \frac{1}{\sqrt
k}(r+is),\frac{w}{\sqrt k}\right),\nonumber
\end{eqnarray}
where $P:\mathbb{C}^\mathrm{n}\times
\mathbb{R}^\mathrm{g}\rightarrow \mathbb{C}$ vanishes to third
order at the origin. Hence, $P_k=:P\left ( \frac{1}{\sqrt
k}(r+is),\frac{w}{\sqrt k}\right)=O(k^{-3/2})$ for fixed $r,s\in
\mathbb{R}^\mathrm{n}$ as $k\rightarrow +\infty$.

On the upshot, we have $u_{k,\varpi}(x+w/\sqrt{k})_{j}\sim
u_{k,\varpi}(x+w/\sqrt{k})_{jb}$ and
\begin{equation}\label{eqn:integral-for-u}
u_{k,\varpi}(x+w/\sqrt{k})_{jb}= \frac{\dim (V_\varpi)
}{(2\pi)^{\mathrm{n}+2}}\,k^{2-\mathrm{g}/2}\,e^{-ik\vartheta_l}\,\int
_{\mathbb{R}^{\mathrm{g}}}\!\int
_{\mathbb{R}^\mathrm{n}}\,F_{jk}(\nu,s)\,d\nu\,ds;
\end{equation}
here, performing the coordinate change $\eta \rightarrow \sqrt k\,
\eta $, we have set
%%%%%%%%%%%%%%%%%%%%%%%%%%%%%%%%
\begin{eqnarray}\label{eqn:Fj(s)}
F_{jk}(\nu,s)&=:&k^{\mathrm{n}/2}\,\int \!\int _0^{+\infty}\!
e^{i\,k\,S}\,A\,dr\,d\theta\, d\tau\,d\eta\, d\vartheta\,dt,
\end{eqnarray}
where the complex phase
%%%%%%%%%%%%%%%%%%%%%%%%%%%%%%%
\begin{eqnarray}\label{phaseSj}
S(\vartheta,t,\theta,\tau,\eta,r)&=:&it[1-e^{i(\vartheta-\theta)}]-\vartheta
+\tau \,\theta+\eta \cdot r\end{eqnarray} satisfies $\Im (S)\ge
0$, and the amplitude $A$ is
%%%%%%%%%%%%%%%%%%%%%%%%%%%%%
\begin{eqnarray}\label{amplj}
A&=:&e^{-i\,t\,e^{i(\vartheta-\theta)}\,\left [2\,\Omega
_{m_j}(\nu _M(m_j),w_j)+iK _2\left (w_j-\nu
_M(m_j),\,r+i\,s\right )\right ]}\nonumber \\
&&e^{i\,k\,t\,e^{i(\vartheta-\theta)}\,P\left ( \frac{1}{\sqrt
k}(r+is),\frac{w}{\sqrt k}\right)}\,\varrho
_{j}(y)\,{f_\lambda}(q)\,D_\Lambda (q)\,\nonumber \\
     && b_{j}(k^{-1/6}\,(\nu,r+is))\,\zeta_{j}\big (\mu _{g^{-1}}\circ
     r_\vartheta
(x_{j}+w_{j}/\sqrt{k}),y,k t\big ) \nonumber \\
     && \gamma _j(e^{i(\vartheta _j+\vartheta)},g_j^{-1}\,e^{\nu/\sqrt k})\,\chi _\varpi
(e^{-\nu /\sqrt k}g_j)\,H_G(\nu/\sqrt k);\end{eqnarray}
%%%%%%%%%%%%%%%%%%%%%%%%%%%%%
in (\ref{eqn:Fj(s)}), the integration from $0$ to $+\infty$ refers
to $dt$. In (\ref{amplj}), $H_G$ denotes the Haar density on $G$,
expressed in the local coordinates given by the exponential chart;
thus, $H_G(0)=1$. The factor $b_{j}(k^{-1/6}\,(\nu,r+is))$ comes
from setting $\xi ={\nu}/{\sqrt k}$ and $z=(r+is)/{\sqrt k}$ in
$b_{jk}$ given by (\ref{eqn:part-fnctn-jk}).

As above, let us write $w_{j}=p_w+iq_w$, with $p_w,q_w\in
\mathbb{R}^\mathrm{n}$. The term $2\,\Omega _{m_j}(\nu
_M(m_j),w_j)+iK _2\left (w_j-\nu _M(m_j),\,r+i\,s\right )$,
appearing in the exponent in the first factor of (\ref{amplj}),
may be rewritten:
\begin{eqnarray}\label{eqn:key-term}
2\, (p_\nu\,q_w-q_\nu\,p_w)-(q_w-q_\nu)\cdot
r+(p_w-p_\nu)\cdot s\nonumber \\
-\frac i2\, (\|p_w-p_\nu-r\|^2+\|q_w-q_\nu-s\|^2).
\end{eqnarray}

\begin{lem}\label{lem:dominated-expansion}
\begin{description}
  \item[i):] As $k\rightarrow +\infty$,
  there is an asymptotic expansion, uniform on compact subsets
  of $\mathbb{R}^\mathrm{g}\times \mathbb{R}^\mathrm{n}$,
  \begin{eqnarray*}F_{jk} (\nu,s)&\sim & k^{\mathrm{n}/2-2}\, Z_{j0}(\nu,s)
  +\sum _{f\ge 1}k^{(\mathrm{n}-f)/2-2}\, Z _{jf} (\nu,s).\end{eqnarray*}
  The coefficient of the leading term is:
\begin{eqnarray*}Z_{j0}(\nu,s)&=&
  \frac{(2\pi)^{\mathrm{n}+2}}{\pi ^\mathrm{n}}\,\chi _\varpi (g_j)\,f_\lambda(x_j)\\
  && e^{-\frac 12\|p_w-p_\nu\|^2}\,
  e^{-i (p_w -p_\nu)s-2\,i\, (p_\nu\,q_w-q_\nu\,p_w)-\frac 12\|q_w-q_\nu
  -s\|^2}\,
  .\end{eqnarray*}
  \item[ii):] there exist positive constants $c>0$ and
  $C_f>0$ for every $f=1,2,\ldots$ such that
  $|Z _{jf}(\nu,s)|\,<\,C_f\,e^{-c(\|s\|^2+\|\nu\|^2)}$, for every
  $(\nu,s)\in \mathbb{R}^\mathrm{g}\times \mathbb{R}^\mathrm{n}$.
  \item[iii):] For every $\ell =0,1,2,\ldots$ there exists $D_\ell>0$ such
  that
  $$\left |F_{jk} (\nu,s)-\sum_{f=0}^\ell k^{(\mathrm{n}-f)/2-2}\, Z
  _{jf}(\nu,s)\right |\le D_\ell \,k^{(\mathrm{n}-\ell-1)/2-2}e^{-c(\|s\|^2+\|\nu\|^2)}$$ for every
  $(\nu,s)\in \mathbb{R}^\mathrm{g}\times \mathbb{R}^\mathrm{n}$.
 \end{description}
 \end{lem}

\textit{Proof.} i): On any fixed compact subset of
$\mathbb{R}^\mathrm{g}\times \mathbb{R}^\mathrm{n}$, we have
$H_G(\nu/\sqrt k)=1+O(k^{-1/2})$, $\mu _{e^{\nu /\sqrt
k}}=\mathrm{id}+O(k^{-1/2})$, $\chi _\varpi
(e^{-\nu/\sqrt{k}}\,g_j)=\chi _\varpi (g_j)+O(k^{-1/2})$ as
$k\rightarrow +\infty$. Incorporating the terms $O(k^{-1/2})$ into
the amplitude, (\ref{eqn:Fj(s)}) may be interpreted as an
oscillatory integral, with complex phase (\ref{phaseSj}), and
whose amplitude may be developed in descending powers of
$k^{-1/2}$.

Since $r=\frac{\partial S}{\partial \eta}$, the asymptotic
contribution to $F_{jk}(s)$ from the region $\|r\|\ge 1$, say, is
$O(k^{-\infty})$.  Therefore, we may assume that both $r$ and $s$
are bounded in norm, and so $b_{j}(k^{-1/6}\,(\nu,v))=1$ if $k\gg
0$. The proof of the following is left to the reader:

\begin{claim} \label{claim:stat-point-sj}
The phase $S$ has only one stationary point $(\vartheta _0
,t_0,\theta _0,\tau_0,\eta_0,r_0)$, given by $t_0=\tau_0=1$,
$\vartheta_0=\theta_0=0$, $r_0=\eta_0=0$. The Hessian of $S$ at
this stationary point is
$$i\left [\begin{array}{cccccc}
1&-i&-1&0&0&0\\
-i&0&i&0&0&0\\
-1&i&1&-i&0&0\\
0&0&-i&0&0&0\\
0&0&0&0&0&-iI_\mathrm{n}\\
0&0&0&0&-iI_\mathrm{n}&0
\end{array}
\right ].$$
\end{claim}

Now (i) follows in view of (\ref{eqn:classical-symb}) and
(\ref{eqn:key-term}) by the complex stationary phase Lemma
(Theorem 7.7.5 of \cite{hor}).

ii): By our choice of adapted Heisenberg local coordinates for
$\Lambda$ centered at $x_{j}$, and by Corollary \ref{cor:xi=0},
$\nu \mapsto p_\nu$ is an injective $\mathbb{R}$-linear map
$\frak{g}\rightarrow \mathbb{R}^\mathrm{n}$; therefore, so is the
affine map $A_j:\frak{g}\oplus \mathbb{R}^\mathrm{n}\rightarrow
\mathbb{C}^\mathrm{n}$ given by
$$A_j(\nu,s)\,=\,: (p_w-p_\nu)\,+\,i(q_w-q_\nu-s).$$ Hence there exist
$c,d>0$ such that
%%%%%%%%%%%%%%%%%%%%%%%%%%%%%%%%%%%%%%%%
\begin{eqnarray}\label{eqn:estimate-xi-nu}
\left \|p_w-p_\nu\right\|^2+\left \|q_w-q_\nu-s\right \|^2\ge
c\left (\|\nu\|^2+\|s\|^2\right )-d.\end{eqnarray}
%%%%%%%%%%%%%%%%%%%%%%%%%%%%%%%%%%%%%%%%%%
Now ii) follows in view of (\ref{eqn:key-term}) and i).

iii): In view of the cut-off $b_{j}(k^{-1/6}\,(\nu,r+is))$, we may
suppose $\|(\nu,s)\|\le k^{1/6}$ (and $\|r\|\le 1$, as above). Let
us make the coordinate change $s'=s/k^{1/6}$, $\nu '=\nu/k^{1/6}$,
so that $\nu'$ and $s'$ may be assumed to be bounded. Then
(\ref{eqn:Fj(s)}) may still be interpreted as an oscillatory
integral, with complex phase $S$, and whose amplitude is
$S_{1/6}$. We may then apply the stationary phase Lemma as in i),
ii) and plug back in $s=k^{1/6}\,s'$ and $\nu =k^{1/6}\,\nu'$ in
the result.

This completes the proof of Lemma \ref{lem:dominated-expansion}.

\bigskip

In view of (\ref{eqn:integral-for-u}), Lemma
\ref{lem:dominated-expansion} implies the asymptotic expansion
\begin{eqnarray}
u_{k,\varpi}(x+w/\sqrt{k})_{jb}&=&k^{(\mathrm{n}-\mathrm{g})/2}\,\varrho
_0
(k,\varpi,w)^{(j)}\nonumber \\
&&+\sum _{f\ge 1}\,k^{(\mathrm{n}-\mathrm{g}-f)/2}\,\varrho _f
(k,\varpi,w,x)^{(j)},
\end{eqnarray}
where
\begin{eqnarray}
\varrho _f (k,\varpi,w,x)^{(j)}&=&\frac{\dim (V_\varpi) }{(2\pi)
^{\mathrm{n}+2} }\,e^{-ik\vartheta
_l}\int_{\mathbb{R}^{\mathrm{g}}}\!\int
_{\mathbb{R}^\mathrm{n}}\,Z_{jf}(\nu,s)\,d\nu\,ds.
\end{eqnarray}
In particular, the coefficient of the leading term is
\begin{eqnarray}\label{eqn:integral-for-leading-term}
\varrho _0 (k,\varpi,w)^{(j)}&=&\frac{\dim (V_\varpi)}{\pi
^\mathrm{n} }\,e^{-ik\vartheta _j}\,\chi _\varpi (g_j)\,f _\lambda
(x_{j})\,\\
  &&\int_{\mathbb{R}^{\mathrm{g}}}\!\int _{\mathbb{R}^\mathrm{n}}\,
  e^{-\frac 12\|p_w-p_\nu\|^2}\,
  e^{-i (p_w-p_\nu )s-2\,i\, (p_\nu\,q_w-q_\nu\,p_w)-\frac 12\|q_w-q_\nu
  -s\|^2}
  \,d\nu\,ds.\nonumber
\end{eqnarray}

To integrate in $ds$, let us make the change of variables
$s\rightsquigarrow s+q_\nu-q_w$, and recall that the function
$e^{-\|x\|^2/2}$ on $\mathbb{R}^\mathrm{n}$ equals its own Fourier
transform. We obtain:
\begin{eqnarray}\label{eqn:integral-for-leading-term-in-ds}
\varrho _0 (k,\varpi,w)^{(j)}&=&\frac{\dim (V_\varpi)}{\pi
^\mathrm{n} }\,(2\pi)^{\frac{\mathrm{n}}{2}}\,e^{-ik\vartheta
_j}\,\chi _\varpi (g_j)\,f _\lambda
(x_{j})\nonumber\\
  &&\int_{\mathbb{R}^{\mathrm{g}}}\,
  e^{-\|p_w-p_\nu\|^2}\,
  e^{-i[(p_w-p_\nu)(q_w-q_\nu)+2\, (p_\nu\,q_w-q_\nu\,p_w)]}
  \,d\nu.
\end{eqnarray}

Let us now decompose $w_j$ as in (\ref{eqn:intundec}), and to
simplify our notation let us write $w_a,\ldots$ for
$(w_j)_a,\ldots$ ($j$ being fixed in our argument). Let us write
$w_a,\ldots$ in local coordinates as column vectors
$(p_{a}\,q_{a})^t,\ldots\in \mathbb{R}^{2n}$. More precisely, we
have:
\begin{claim}\label{claim:wj}
%%%%%%%%%%%%%%%%%%%%%%%%%
%%%%%%%%%%%%%%%%%%%%%%%
\begin{description}
  \item[i):] $w_a=\begin{pmatrix}
  p_a \\
  0
\end{pmatrix}$ where $p_{a}\cdot p_\nu=0\,\,\forall
\,\nu \in \frak{g}$;
  \item[ii):] $w_b=\begin{pmatrix}
  0 \\
  q_{b}
\end{pmatrix}$ and $w_c=\begin{pmatrix}
  0 \\
  q_{c}
\end{pmatrix}$  for appropriate $q_{b},\,q_c \in \mathbb{R}^n$,
and $p_\nu\cdot q_{c}=0$, $\forall\,\nu \in \frak{g}$;
  \item[iii):]$ w_d=\begin{pmatrix}
  p_{\nu_w} \\
  q_{\nu_w}
\end{pmatrix},$ for a unique
$\nu_w \in \frak{g}$.
\end{description}
\end{claim}

\textit{Proof.} For the second statement of ii), recall that
$w_c\in T_{m_j}\Lambda '\subseteq T_{m_j}M'$, and the latter is
the symplectic annihilator of $\frak{g}_M(m_j)$. Everything else
is immediate.

\begin{lem}\label{lem:symplecticpairing}
Let $p_{w_j},\,q_{w_j}\in \mathbb{R}^n$ be as in Claim
\ref{claim:wj}. Then for every $\nu \in \frak{g}$ one has
\begin{equation}\label{eqn:symppairing}
p_\nu\cdot q_w-q_\nu\cdot p_w=p_\nu \cdot q_{b}-q_\nu\cdot
p_{a}.\end{equation}
\end{lem}

\textit{Proof.} The left hand side of (\ref{eqn:symppairing}) is
the symplectic pairing $\Omega _{m_j}\big (\nu _M(m_j),w_j\big )$.
We have $p_w=p_{a}+p_{\nu _w}$, $q_w=q_{b}+q_{c}+q_{\nu _w}$.
However, recalling that $\frak{g}_M(m_j)\subseteq T_{m_j}M$ is an
isotropic subspace, we have $\Omega _{m_j}\big (\nu _M(m_j),\nu
_w(m_j)\big )=0$ for every $\nu \in \frak{g}$; therefore, $p_{\nu
_w}$ and $q_{\nu _w}$ may be ignored. The statement then follows
from ii) of Claim \ref{claim:wj}.

\bigskip
Let us make the change of variables $\beta =\nu-\nu_w$ in
(\ref{eqn:integral-for-leading-term-in-ds}). The real part of the
exponent in (\ref{eqn:integral-for-leading-term-in-ds}) then is
$-\|p_a\|^2-\|p_\beta\|^2$, while the imaginary part may be
written as
\begin{eqnarray}\label{eqn:exponent-in-integral-for-leading-term-in-ds}
-\big (p_a-p_\beta\big )\,\cdot \big[(q_b+q_c)-q_\beta\big] -2\big
[(p_\beta+p_{\nu _w})\cdot q_b-(q_\beta+q_{\nu _w})\cdot p_a\big
]\nonumber\\
=\left [\Omega _{m_j}(w_j-w_d,w_a)+2\,\Omega_{m_j}(w_j,w_d)\right
] -\big (p_\beta\cdot q_b-3q_\beta \cdot p_a\big ) -p_\beta\cdot
q_\beta .
\end{eqnarray}

Let us set $C(w_j)=\Omega
_{m_j}(w_j-w_d,w_a)+2\,\Omega_{m_j}(w_j,w_d)$. We may then rewrite
the right hand side of (\ref{eqn:integral-for-leading-term-in-ds})
as
\begin{eqnarray}\label{eqn:new-integral-for-leading-term-in-ds}
\varrho _0 (k,\varpi,w)^{(j)}&=&\dim (V_\varpi)\,(2\pi)^{\frac 32
\mathrm{n}+1}\,e^{-ik\vartheta _j}\,\chi _\varpi (g_j)\,f _\lambda
(x_{j})\,e^{-\|w_a\|^2+iC(w_j)}\nonumber\\
  &&\int_{\mathbb{R}^{\mathrm{g}}}\,
  e^{-\big (\|p_\beta\|^2+i\,p_\beta\cdot q_\beta \big )}\,
  e^{-i\big (p_\beta\cdot q_b-3q_\beta \cdot p_a\big )}
  \,d\beta.
\end{eqnarray}
In (\ref{eqn:new-integral-for-leading-term-in-ds}), the integral
is over $\frak{g}$, identified with $\mathbb{R}^{\mathrm{g}}$ by
means of an orthonormal basis for the Haar metric. Thus, the
Lebesgue measure $d\beta$ corresponds to the Haar measure at the
identity $e\in G$. To make our statement more intrinsic, we shall
now rewrite the latter integral as an integral over $\frak{g}
_M(m_j)\subseteq T_{m_j}M$, with the induced metric.

\begin{lem}\label{lem:effective-volume}
Fix $t\in M'$. Suppose that $\mathcal{B}$ is an orthonormal basis
of $\frak{g}$ for the Haar metric, and that $\mathcal{B}_t$ is an
orthonormal basis of $\frak{g} _M(t)$ for the induced metric from
$T_tM$. Identify $\frak{g}$ with $\frak{g} _M(t)$ by the linear
isomorphism $\xi \mapsto \xi _M(t)$. Let
$A=M^{\mathcal{B}_t}_{\mathcal{B}}(\mathrm{id}_{\frak{g}})$ be the
matrix of the base change. Then $$|\det
(A)|=\frac{1}{V_{\mathrm{eff}}(t)\,|G_{t}|},$$ where
$V_{\mathrm{eff}}(t)$ is the effective potential at $t$, and
$G_{t}\subseteq G$ is the stabilizer subgroup of $t$.
\end{lem}

\begin{rem} Since $m_j=\mu _{g_j}(m)$, we have
$V_{\mathrm{eff}}(m_j)=V_{\mathrm{eff}}(m)$ and
$|G_{m_j}|=|G_m|$.\end{rem}

\textit{Proof.} Suppose $\mathcal{B}=\{v_1,\ldots,v_\mathrm{g}\}$,
$\mathcal{B}_t=\{w_1,\ldots,w_\mathrm{g}\}$ so that $w_j=\sum
_{i=1}^\mathrm{g}a_{ij}v_i$, where $A=[a_{ij}]$. Hence
\begin{equation}\label{eqn:det-and-change} w_1\wedge \cdots \wedge
w_\mathrm{g}=\det (A)v_1\wedge \cdots \wedge
v_\mathrm{g}.\end{equation} Let $\mathrm{dens}_G$ be the Haar
density on $G$, so that $\int _G\mathrm{dens}_G=1$; hence
$\mathrm{dens}_G(v_1\wedge \cdots \wedge v_\mathrm{g})=1$. Let
$\mathrm{dens}_t$ be the pull-back to $G$ of the invariant metric
density on the orbit $G\cdot t\cong G/G_t$ under the degree -
$|G_t|$ covering map $g\mapsto g\cdot t$. By invariance,
\begin{equation}\label{eqn:veff-and-change}
\mathrm{dens}_t=V_{\mathrm{eff}}(t)\cdot |G_t|\cdot
\mathrm{dens}_G.\end{equation} By construction,
\begin{eqnarray*}1=\mathrm{dens}_t(w_1\wedge \cdots \wedge
w_\mathrm{g})&=&|\det
(A)|\cdot V_{\mathrm{eff}}(t)\cdot |G_t|\cdot \mathrm{vol}_G(v_1\wedge \cdots \wedge v_\mathrm{g})\\
&=&|\det (A)|\cdot V_{\mathrm{eff}}(t)\cdot |G_t|.\end{eqnarray*}

\bigskip
Now let $\beta$ and $u$ denote the linear coordinates on
$\frak{g}$ associated to the basis $\mathcal{B}$ and
$\mathcal{B}_q$, respectively; thus, $\beta =A\,u$. By the Lemma,
\begin{equation}\label{eqn:change-of-linear-coordinates}
d\beta =\big (V_{\mathrm{eff}}(q)\,|G_q|\big )
^{-1}\,du.\end{equation}

We have already exploited the following consequence of Corollary
\ref{cor:xi=0}: working in adapted Heisenberg local coordinates,
the projection of
$$\frak{g}_M(m_j) \subseteq T_{m_j}M\cong \mathbb{C}^\mathrm{n}\cong
\mathbb{R}^\mathrm{n}\times \mathbb{R}^\mathrm{n}$$ onto
$\mathbb{R}^\mathrm{n}\times \{0\}$,
$$\mathrm{p}_{\mathrm{real}}:\nu _M(m_j)\mapsto p_\nu,\,\,\,\,\,\,(\nu\in \frak{g})$$
is injective. Hence there exists a linear map
$T:\mathbb{R}^\mathrm{n}\rightarrow \mathbb{R}^\mathrm{n}$ such
that $q_\nu=T(p_\nu)$ for every $\nu \in \frak{g}$ (if we so wish,
we may determine $T$ uniquely by imposing that it vanishes on the
Euclidean orthocomplement of
$\mathrm{p}_{\mathrm{real}}(\frak{g}_M(p_j) )\subseteq
\mathbb{R}^n$). We shall think of $T$ as an n$\times$n real
matrix.

On the upshot, identifying $\frak{g}\cong \mathbb{R}^\mathrm{g}$
by $ \mathcal{B}_{m_j}$, and $T_{m_j}M\cong
\mathbb{C}^{\mathrm{n}}$ by the given choice of adapted Heisenberg
local coordinates, the inclusion $ \frak{g}\cong
\frak{g}_{M}(m_j)\hookrightarrow T_{m_j}M$ may be written
%%%%%%%%%%%%%%%
%%%%%%%%%%%%%%%
\begin{equation}\label{eqn:RandT}
\iota (u)\,=\,Ru+iTRu\,\,\,\,\,\,(u\in \mathbb{R}^{\mathrm{g}}),
\end{equation}
%%%%%%%%%%%%%%%
%%%%%%%%%%%%%%%
for a certain $\mathrm{n}\times \mathrm{g}$ real matrix $R$ of
maximal rank g. Since $ \mathcal{B}_{m_j}$ is orthonormal for the
induced metric, we have $RR^t+R^tT^tTR=\mathrm{I}_{\mathrm{g}}$.
%%%%%%%%%%%%%%%%%%%%%%%%%%%%%%
%%%%%%%%%%%%%%%%%%%%%%%%%%%%%
\begin{lem}\label{lem:symmetric-matrix}
$R^tTR$ is a $\mathrm{g}\times \mathrm{g}$ symmetric
matrix.\end{lem}

\textit{Proof.} Since $m_j\in M'$, $\frak{g}_M(m_j)\subseteq
T_{m_j}M$ is an isotropic subspace. Thus, for every $\xi,\,\nu\in
\frak{g}$ we have \begin{eqnarray*}\Omega _{m_j}\big (\xi
_M(m_j),\nu_M(m_j)\big )=\Omega _0\left (\begin{pmatrix}
  p_\xi \\
  q_\xi
\end{pmatrix},\begin{pmatrix}
  p_\nu \\
  q_\nu
\end{pmatrix}\right )=\Omega _0\left (\begin{pmatrix}
  R\xi \\
  TR\xi
\end{pmatrix},\begin{pmatrix}
  R\nu \\
  TR\nu
\end{pmatrix}\right )\\
=\xi ^t(R^tTR-R^tT^tR)\nu,\end{eqnarray*} where $\Omega_0$ denotes
the standard symplectic structure on $
\mathbb{R}^{\mathrm{n}}\times \mathbb{R}^{\mathrm{n}}$.

\bigskip

Now (\ref{eqn:new-integral-for-leading-term-in-ds}) may be
rewritten as
%%%%%%%%%%%%%%%%%%%%%%%%%%%%
\begin{eqnarray}\label{eqn:new-integral-for-leading-term-in-ds-in-du}
\varrho _0 (k,\varpi,w)^{(j)}&=&\frac{\dim
(V_\varpi)}{V_{\mathrm{eff}}(q)\,|G_q|}\,(2\pi)^{\frac 32
\mathrm{n}+1}\,e^{-ik\vartheta _j}\,\chi _\varpi (g_j)\,f _\lambda
(x_{j})\,e^{-\|w_a\|^2+iC(w_j)}\nonumber\\
  &&\int_{\mathbb{R}^{\mathrm{g}}}\,
  e^{-u^t\big (R^tR+i\,R^t T R\big )u}\,
  e^{-iu^tR^t\big ( q_b-3T^t p_a\big )}
  \,du.
\end{eqnarray}

Up to a scalar factor $(2\pi) ^{\frac {\mathrm{g}}{2}}$, the
integral in (\ref{eqn:new-integral-for-leading-term-in-ds-in-du})
is the evaluation at $R^t (q_b-3T^t p_a)$ of the Fourier transform
of the function $e^{-\frac 12(u,Au)}$ on $ \mathbb{R}^\mathrm{g}$,
where $A=2(R^tR+i\, R^tTR)$ is a $\mathrm{g}\times \mathrm{g}$
complex symmetric matrix with positive definite real part. By
Theorem 7.6.1 of \cite{hor}, we have

\begin{eqnarray}\label{eqn:new-integral-for-leading-term-fourier-trasform}
\varrho _0 (k,\varpi,w)^{(j)}&=&\frac{\dim
(V_\varpi)}{V_{\mathrm{eff}}(q)\,|G_q|}\,\frac{1}{\pi ^\mathrm{n}
}\,\sqrt{\frac{(2\pi)^{
 (\mathrm{n}+\mathrm{g})}}{2^{g}}}\, e^{-ik\vartheta _j}\,\chi
_\varpi (g_j)\,f _\lambda
(x_{j})\nonumber \\
  &&\det \big (R^tR+i\,R^t T R\big )^{-1/2}\,\exp (-Q(w_j))\end{eqnarray}
%%%%%%%%%%%%%%
%%%%%%%%%%%%%%%
where

  \begin{eqnarray}\label{eqn:local-normal-gauge}
  Q(w_j)&=&\|p_a\|^2-i C(w_j)+\frac 14 \left ( R^t (q_b-3T^t p_a),
  \big (R^tR+i\,R^t T R\big )^{-1}R^t (q_b-3T^t
  p_a)\right)\nonumber\\
  &=&S(w'_j)+i P(w_j),
\end{eqnarray}
where $S$ and $P$ denote real valued quadratic forms. Here
$(\,,\,)$ is the standard Euclidean scalar product on
  $\mathbb{R}^\mathrm{g}$. Thus, if
$\big (R^tR+i\,R^t T R\big )^{-1}=F+iG$, where $F$ and $G$ are
g$\times$g real symmetric matrices,
\begin{equation}\label{eqn:transv-quadr-form}
S(w'_j)=:\|p_a\|^2+\frac 14\, \left ( R^t (q_b-3T^t p_a),
  F\,R^t (q_b-3T^t
  p_a)\right).\end{equation}
    Recall that $w'_j=w_a+w_b$ in the decomposition described in
  Definition \ref{defn:general-defn} and Claim \ref{claim:wj}.

  \begin{lem} $S(w_j')\ge 0$, and equality holds only if $w_j'=0$.
  \label{lem:non-deg-quad}
  \end{lem}

\textit{Proof.} Since $F$ is positive definite by construction,
both summands in (\ref{eqn:transv-quadr-form}) are $\ge 0$.
Suppose $S(w_j')=0$. Then both summands vanish, whence $p_a=0$ and
$R^tq_b=0$. Thus we are reduced to proving:

\begin{lem}\label{lem:non-deg-quad-sub}
If $R^tq_b=0$, then $q_b=0$.
\end{lem}

\textit{Proof.} By construction, the range of $R$ is
$$p_{\mathrm{real}}\left (\frak{g}_M(m_j)\right )=\{p_\nu :\nu \in
\frak{g}\}\subseteq \mathbb{R}^\mathrm{n}.$$ Recall that the
symplectic annihilator of $\frak{g}_M(m_j)$ is given by
$\frak{g}_M(m_j)^0=T_{m_j}M'$. Hence, in view of the
identification $T_{m_j}M\cong \mathbb{R}^\mathrm{n}\oplus
\mathbb{R}^\mathrm{n}$ (and viewing as usual $T_{x_j}\Lambda$ as a
subspace of $T_{m_j}M$),
\begin{eqnarray*}\ker (R^t)&=&\left \{q\in
\mathbb{R}^\mathrm{n}:p_\nu^t\,q=0\, \forall \nu \in
\frak{g}\right \}=\left \{q\in
\mathbb{R}^\mathrm{n}:\begin{pmatrix}
  0 \\
  q
\end{pmatrix}\in \frak{g}_M(m_j)^{0 }\right \}\\
&=&\left \{q\in \mathbb{R}^\mathrm{n}:\begin{pmatrix}
  0 \\
  q
\end{pmatrix}\in T_{m_j}M'\right \}=\left \{q\in \mathbb{R}^\mathrm{n}:\begin{pmatrix}
  0 \\
  q
\end{pmatrix}\in T_{x_j}\Lambda'\right \}.\end{eqnarray*}
%%%%%%%%%%%%%%%%%%
By definition, $w_b\in (T_{x_j}\Lambda ')^\perp$. Thus if
$R^tq_b=0$ then
$$
\begin{pmatrix}
  0 \\
  q_b
\end{pmatrix}\in T_{x_j}\Lambda'\cap (T_{x_j}\Lambda
')^\perp=\{0\}.$$

\bigskip

This completes the proof of Theorem \ref{thm:main2}.

\begin{rem}\label{rem:negative-answer}
Let us now consider the question described in the introduction,
i.e. whether quantization commutes with reduction for a transverse
compact Legendrian submanifold $\Lambda \subseteq X$. For the sake
of brevity, we shall make fairly brutal simplifying assumptions,
and leave it to the interested reader to work out more general
cases.

Let us assume that $G$ acts freely on $M'$, and - to fix ideas -
that $\Lambda$ meets every $S^1\times G$-orbit in $X'$ at most
once. Thus, $\Lambda '=:\Lambda \cap X'$ is an
$(\mathrm{n-g})$-dimensional isotropic submanifold, which projects
diffeomorphically onto a compact Legendrian submanifold $\Lambda
_0\subseteq X_0=:X'/G$ ($\Lambda$ and $\Lambda_0$ actually map
down diffeomorphically onto Lagrangian submanifolds in $M$ and
$M_0$, respectively). Also, let us choose as a half-density on
$\Lambda$ the Riemannian half-density
$\mathrm{dens}_\Lambda^{(1/2)}$, so that $f_\lambda=1$.

Let us fix $x\in \Lambda '$, and let $x_0\in X_0$ be its image in
$X_0$. To simplify, let us also assume that $\Lambda$ is
perpendicular to the $G$ orbit $G\cdot x$ at $x$, so that
-referring to (\ref{eqn:RandT}) - we have $T=0$ and
$R^tR=I_\mathrm{g}$. In view of Theorem \ref{thm:main2}, we then
have:
%%%%%%%%%%%%%%%%%%%%%%%%%%%%%%%%%%%%%%%%%%%
%%%%%%%%%%%%%%%%%%%%%%%%%%%%%%%%%%%%%%%%%%%
\begin{eqnarray*}
  u_{k,0}(x)&\sim &
  \frac{\pi ^{-\mathrm{n}}}{V_{\mathrm{eff}}\big (\pi (x)\big )}\,
  \sqrt{\frac{(2\pi)^{\mathrm{n}+\mathrm{g}}}{2^{\mathrm{g}}}}
  \,k^{(\mathrm{n-g})/2}
  +\sum _{f\ge
  1}\varrho_f\,k^{(\mathrm{n-g}-f)/2}.
  \end{eqnarray*}
%%%%%%%%%%%%%%%%%%%%%%%%%%%%%%%%%%%%%%%%%%%%%
%%%%%%%%%%%%%%%%%%%%%%%%%%%%%%%%%%%%%%%%%%%%%

Now there are two natural ways to induce a half-density on
$\Lambda'\cong\Lambda_0$: One is to choose the Riemannian
half-density, $\lambda'= \mathrm{dens}_{\Lambda_0}^{(1/2)}$, so
that $f_{\lambda'}=1$. The other is to divide $
\mathrm{dens}_\Lambda^{(1/2)}$ by the half-density on
$\frak{g}^*\cong \frak{g}$ associated to the Haar metric. Let
$\lambda''$ be the half-density obtained in this manner. By
arguments similar to those in Lemma \ref{lem:effective-volume},
one can see that $f_{\lambda''}(x_0)=V_\mathrm{eff}\big
(\pi(x)\big )^{-1/2}$.

If $u_k'\in H^0(M_0,L_0^{\otimes k})\cong \mathcal{H}_k(X_0)$ is
the sequence associated to $\lambda'$, then by Corollary
\ref{cor:action-free-case} we have:
%%%%%%%%%%%%%%%%%%%%%%%%%%%%%
%%%%%%%%%%%%%%%%%%%%%%%%%%%%%%%
%%%%%%%%%%%%%%%%%%%%%%%%%%%%%%%
\begin{eqnarray*}
  u_{k}'(x)&\sim &
  \frac{(2\pi)^{\frac{\mathrm{n}}{2}}}{\pi ^\mathrm{n}}\,k^{(\mathrm{n-g})/2}+\sum _{f\ge
  1}\varrho_f\,k^{(\mathrm{n}-f)/2}.
  \end{eqnarray*}
If on the other hand $u_k''\in H^0(M_0,L_0^{\otimes k})\cong
\mathcal{H}_k(X_0)$ is the sequence associated to $\lambda''$, the
leading order term gets multiplied by $V_\mathrm{eff}\big
(\pi(x)\big )^{-1/2}$.
\end{rem}

\section{The Hermitian products.}\label{sctn:herm-prdcts}

%\subsection{Statement of the problem and preliminaries.}
Let us now assume that $\Lambda,\,\Sigma\subseteq X$ are two
compact Legendrian submanifolds, and that $\lambda$ and $\sigma$
are given smooth half-densities on $\Lambda$ and $\Sigma$,
respectively. Let $u=:\Pi _X(\delta _{\Lambda,\lambda}),\,v=:\Pi
_X(\delta _{\Sigma,\sigma})$. Let as usual $\varpi$ be a fixed
highest weight of $G$. We shall study in this section the
asymptotics of the Hermitian products
$(u_{k,\varpi},v_{k,\varpi})_{L^2(X)}$ as $k\rightarrow +\infty$.

In the action-free case, we shall reproduce expansions similar to
those in \cite{bpu}, except for some differences due to the fact
that we are dealing with half-densities rather than half-forms.

To this end, let us first of all recall that we are unitarily and
equivariantly identifying functions and half-densities on $X$.
Furthermore, the self-duality pairing $<\,,\,>$ and the
$L^2$-unitary product $(\,,\,)_{L^2}$ of smooth half-densities
$\tau=f\cdot \mathrm{dens}_X$ and $\upsilon =g\cdot
\mathrm{dens}_X$ are related by $\int _Xf\cdot
\overline{g}\,\mathrm{dens}_X=\left (\tau,\upsilon\right
)_{L^2}=\left <\tau,\overline{\upsilon}\right
>$. Suppose then that $u_t$, $t>0 $, is a family of smooth
half-densities on $X$ such that $u_t\rightarrow u$ as
$t\rightarrow 0$, in the topology of the space of all generalized
half-densities whose wave front is conormal to $\Lambda$. In view
of the self-adjointness of the orthogonal projector $\Pi
_{k,\varpi}$ on $\mathcal{H}(X)_{k,\varpi}$, we obtain:

\begin{eqnarray}\label{eqn:reproducing}
(u_{k,\varpi},v_{k,\varpi})_{L^2(X)}&=&\lim _{t\rightarrow 0}\Big
(\Pi _{k,\varpi}(u_t),v_{k,\varpi})\Big )_{L^2(X)}=\lim
_{t\rightarrow
0}\big(u_t,v_{k,\varpi}\big)_{L^2(X)}\nonumber \\
&=&\lim _{t\rightarrow 0} \big<u_t,\overline{v_{k,\varpi}}
\big>=\big <u,\overline{v_{k,\varpi}}
\big>\nonumber \\
&=&\int _\Lambda\, f_\lambda \cdot
\overline{v_{k,\varpi}}\,\mathrm{dens} _\Lambda,\end{eqnarray}
where $\mathrm{dens} _\Lambda$ is the Riemannian density on
$\Lambda$.

\subsection{The transverse case.}\label{sctn:transversecase}
Consider the smooth map given by group action restricted to
$\Lambda$,
$$\Upsilon:
(h,g,x)\in S^1\times G\times \Lambda\mapsto (h,g)\cdot x\in X,$$
To fix ideas, suppose first that $\Upsilon$ is transversal to
$\Sigma'=\Sigma \cap \Big((\Phi \circ \pi )^{-1}(0)\Big)$. In this
case, $\Upsilon^{-1}(\Sigma')$ is a finite set:
$$\Upsilon^{-1}(\Sigma')=\{\tilde y_1,\ldots,\tilde y_r\},$$ where
$\tilde y_j=(h_j,g_j,y_j)$ for some $h_j\in S^1$, $g_j\in G$ and
$y_j\in \Lambda'$. Hence $\widehat{y_j}=:\Upsilon ( \tilde
y_j)=(h_j,g_j)\cdot \tilde y_j\in \Sigma'$ for every $j$.

Now let $U_j\subseteq \Lambda$ be some arbitrarily small
neighbourhood of $y_j$. Since $v_{k,\varpi}=O(k^{-\infty})$ away
from $\Sigma'$, in view of (\ref{eqn:reproducing}) we have
\begin{equation}\label{eqn:reproducing-asympt-transv}
(u_{k,\varpi},v_{k,\varpi})_{L^2(X)}\sim \sum _{j=1}^r\int
_{U_j}\, f_\lambda \cdot
\overline{v_{k,\varpi}}\,\mathrm{dens}_\Lambda.\end{equation}

Let us fix Heisenberg local coordinates $(p,q,\theta)$ for $X$
centered at $\widehat{y}_j$ and adapted to $\Sigma$, defined on an
open neighbourhood $V_j\ni \widehat{y}_j$. Thus, $\Sigma \cap
V_j\subseteq V_j$ is defined by conditions $\theta =f(q)$ and
$p=h(q)$, as described in \S \ref{subsctn:adapted}. We may
arrange, given our assumptions, that
\begin{equation}\label{eqn:arrange-1}
T_{\widehat{y}_j}\Sigma '=\mathrm{span}\left \{
\left.\frac{\partial}{\partial q_1}\right
|_{\widehat{y}_j},\cdots,\left.\frac{\partial}{\partial
q_{\mathrm{n-g}}}\right |_{\widehat{y}_j}\right \}.
\end{equation}

%%%%%%%%%%%%%%%%%%%
%%%%%%%%%%%%%%%%%%%
%%%%%%%%%%%%%%%%%%%
The following is left to the reader:
\begin{lem}\label{lem:tang-space-to-orbit}
Given (\ref{eqn:arrange-1}), we have
\begin{equation}\label{eqn:arrange-2}
\frak{g}_X(\widehat{y}_j) =\mathrm{span}\left \{
\left.\frac{\partial}{\partial p_{\mathrm{n-g}+1}}\right
|_{\widehat{y}_j}+t_{\mathrm{n-g}+1},\cdots,\left.\frac{\partial}{\partial
p_{\mathrm{n}}}\right |_{\widehat{y}_j}+t_{\mathrm{n}}\right \},
\end{equation}
for appropriate $t_{\mathrm{n-g}+1},\ldots,t_\mathrm{n}\in
T_{\widehat{y}_j}\Sigma '$.
\end{lem}

Let us now consider the Legendrian submanifold
\begin{equation*}
\widehat{y}_j\in\Lambda _j=:\Upsilon \big (\{(h_j,g_j)\}\times
\Lambda\big )\subseteq X,\end{equation*} obtained by 'translating'
$\Lambda$ by the action of $(h_j,g_j)\in S^1\times G$. Given
(\ref{eqn:arrange-1}) and Lemma \ref{lem:tang-space-to-orbit}, the
present transversality assumption implies:

\begin{lem}\label{lem:local-coordinates}
In the above situation, $(p_1,\ldots,p_{\mathrm{n-g}})$ restrict
to local coordinates on $\Lambda'_j$ centered at $\widehat{y}_j$,
and
$(p_1,\ldots,p_{\mathrm{n-g}},q_{\mathrm{n-g}+1},\ldots,q_\mathrm{n})$
restrict to local coordinates on $\Lambda _j$ centered at
$\widehat{y}_j$.
\end{lem}

Therefore
$(p_1,\ldots,p_{\mathrm{n-g}},q_{\mathrm{n-g}+1},\ldots,q_\mathrm{n})$
may be viewed in a natural manner as local coordinates on
$\Lambda$ centered at $y_j$, defined on some open neighbourhood
$U_j\subseteq \Lambda$. In order to apply Theorem \ref{thm:main2},
we need to relate these coordinates on $\Lambda$ to the local
Heisenberg coordinates on $X$. Given
$x=(x_1,\ldots,x_\mathrm{n})$, to simplify our notation let us
write $x=(x',x'')$, where $x'=(x_1,\ldots,x_{\mathrm{n-g}})$,
$x''=(x_{\mathrm{n-g+1}},\ldots,x_{\mathrm{n}})$. The following is
left to the reader:

We have:

\begin{lem}\label{lem:new-coordinates}
There exists an $\mathbb{R}$-linear map
$$A_j:\mathbb{R}^\mathrm{n}\rightarrow \mathbb{C}^\mathrm{n}\cong
\{0\}\oplus\mathbb{C}^\mathrm{n} \subseteq \mathbb{R}\oplus
\mathbb{C}^\mathrm{n}$$ such that if $y\in U_j\subseteq \Lambda$
has local coordinates $\frac{1}{\sqrt k}(p',q'')$ on $\Lambda$,
then it has local Heisenberg coordinates $\frac{1}{\sqrt
k}A_j(p',q'')+O(k^{-1})$.
\end{lem}

Let $y\left ( \frac{1}{\sqrt k}(p',q'')\right )$ denote the point
in $U_j$ having local coordinates $\frac{1}{\sqrt k}(p',q'')$. By
Theorem \ref{thm:main2} and Lemma \ref{lem:new-coordinates},
passing to rescaled coordinates on $U_j$ we may then write the
$j$-th summand in (\ref{eqn:reproducing-asympt-transv}) as:

\begin{eqnarray}\label{eqn:rescaled-local-contr}
\int _{U_j}\, f_\lambda \cdot
\overline{v_{k,\varpi}}\,\mathrm{dens}_\Lambda
&=&k^{-\mathrm{n}/2}\int _{\mathbb{R}^\mathrm{n}}f_\lambda \left (
(k^{-1/2}(p',q'')\right ) \,\overline{v_{k,\varpi}}\left
(k^{-1/2}\,A_j(p',q'')+O(k^{-1})\right )\nonumber\\
&&D_\Lambda \big (k^{-1/2}(p',q'')\big )\,dp'\,dq''
\end{eqnarray}
Inserting the asymptotic expansion of Theorem \ref{thm:main2} in
(\ref{eqn:rescaled-local-contr}), we conclude

%%%%%%%%%%%%%%%%%%%%%%%%%%%%%%%%%
%%%%%%%%%%%%%%%%%%%%%%%%%%%%%%%%%%

\begin{prop}\label{prop:transverse-case}
If $\Upsilon :S^1\times G\times \Lambda\rightarrow X$ is
transversal to $\Sigma'$, the $j$-th summand in
(\ref{eqn:reproducing-asympt}) is
\begin{eqnarray*}
\int _{U_j}\, f_\lambda \cdot
\overline{v_{k,\varpi}}\,\mathrm{dens}_\Lambda
&\sim&k^{-\mathrm{g}/2}\rho _0^{(j)}+\sum _{f\ge
  1}k^{-(\mathrm{g}+f)/2}\,\rho_f^{(j)},
\end{eqnarray*}
where
\begin{eqnarray*}
\rho _0^{(j)}&=&
 \frac{\dim(V_\varpi)}{|G_{\pi (y_j)}|}\,
 \frac{1}{\pi ^{\mathrm{n}}}\,\sqrt{\frac{(2\pi)^{\mathrm{n}+\mathrm{g}}}{2^{\mathrm{g}}}}
  \,h_{j}^{k}\,\chi _\varpi (g_j^{-1})\,\overline{\Xi _{\Lambda}
  (\widehat{y_j})}\,
  f_\lambda (y_j)\,\overline{f_\sigma(\widehat{y_j}) }\nonumber\\
  &&\,\cdot \int _{\mathbb{R}^{\mathrm{n}}}
  e^{-S_{\widehat{y_j}}\big (A_j(p',q''),A_j(p',q'')\big )-
  i\,T_{\widehat{y_j}}\big(A_j(p',q''),A_j(p',q'')\big)}\,dp'\,dq''.
\end{eqnarray*}
\end{prop}

%%%%%%%%%%%%%%%%%%%%%%%%%%%%%%%%%
%%%%%%%%%%%%%%%%%%%%%%%%%%%%%%%%%%

In the action-free case, the present transversality assumption
means that $S^1\times \Lambda\rightarrow X$ is transverse to
$\Sigma$. For every $j=1,\ldots,r$, $T_{\pi
(\widehat{y}_j)}\Lambda _j\subseteq T_{\pi (\widehat{y}_j)}M$ is a
Lagrangian subspace transversal to $T_{\pi
(\widehat{y}_j)}\Sigma$. Thus, in the given Heisenberg local
coordinates adapted to $\Sigma$ at $\widehat{y}_j$, we have
\begin{equation}\label{eqn:tag-space-and-symm-matrix}
T_{\pi (\widehat{y}_j)}\Lambda _j=\{(p,Z_jp):p\in
\mathbb{R}^n\}\subseteq T_{\pi (\widehat{y}_j)}M\cong
\mathbb{R}^\mathrm{n}\oplus \mathbb{R}^\mathrm{n},\end{equation}
where $Z_j$ is a symmetric matrix. Therefore, the $p$'s restrict
to a system of local coordinates on $\Lambda_j$ (whence on
$\Lambda$), and $A_j(p)=p+iZ_jp$.

Let $\imath
_{J_{\pi(\widehat{y}_j)}}:\mathrm{Gr}_{\mathrm{lag}}(T_{\pi(\widehat{y}_j)}M)\times
\mathrm{Gr}_{\mathrm{lag}}(T_{\pi(\widehat{y}_j)}M)\rightarrow
\mathbb{R}$ be the invariant introduced in section
\ref{subsctn:lag-invariant}; let us write
$J_j=J_{\pi(\widehat{y}_j)}$. Applying the asymptotic expansion of
Corollary \ref{cor:action-free-case},

\begin{cor}\label{cor:paring-transv-action-free}
Suppose that the two projections $\Lambda\rightarrow M$ and
$\Sigma\rightarrow M$ are transversal. Let $ \Upsilon:S^1\times
\Lambda \rightarrow \Sigma$ be the map induced by the action, and
suppose $\Upsilon^{-1}(\Sigma)=\{\tilde y_1,\ldots,\tilde y_r \}$,
where $\tilde y_j=(h_j,y_j)$. Set $\widehat{y}_j=:h_j\cdot y_j$
and $\Lambda _j=:h_j\cdot \Lambda$ for every $j$. Then
\begin{equation*}
(u_k,v_k)
 \sim \rho _0+\sum _{f\ge
  1}k^{-f/2}\,\rho_f,
\end{equation*}
where \begin{eqnarray*}\rho _0&=&
  \frac{ (2\pi) ^{\frac{\mathrm{n}}{2}}   }{      \pi ^{\mathrm{n}}       }\,
  \sum _{j=1}^r\,h_{j}^{k}\,
  f_\lambda (y_j)\,\overline{f_\sigma(\widehat{y_j}) }\,
  \imath _{J_j}\big (T_{\widehat{y}_j}\Lambda_j,
  T_{\widehat{y}_j}\Sigma\big ) ^{-1}\int _{\mathbb{R}^{\mathrm{n}}}
  e^{-\|p\|^2+
  i\,p^tZ_jp}\,dp.
\end{eqnarray*}
\end{cor}

\subsection{The clean case.}

Now we shall make the following more general hypothesis:

\begin{description}
    \item[i):] $\Lambda$ and $\Sigma$ are both transversal to
    $X'$; let us set $\Lambda '=:\Lambda\cap X'$,
    $\Sigma'=:\Sigma\cap X'$.
    \item[ii):] the smooth map given by group action restricted to $\Lambda$,
$$\Upsilon:
(h,g,x)\in S^1\times G\times \Lambda\mapsto (h,g)\cdot x\in X,$$
meets $\Sigma'$ nicely; by this, we mean that every connected
component of $\Upsilon ^{-1}(\Sigma')$ is a manifold, and that for
every $\varsigma =(h,g,x)\in \Upsilon ^{-1}(\Sigma')$ we have
$$T_\varsigma\left (\Upsilon ^{-1}(\Sigma')\right
)=(d_\varsigma\Upsilon)^{-1}\left (
T_{\Upsilon(\varsigma)}\Sigma'\right ).$$

    \item[iii):] there exist integers $r_\Sigma, r_\Lambda\ge 1$
    such that for every $x\in \Sigma '$ and $y\in \Lambda'$ one has
    \begin{equation}\label{eqn:constant-cardinality}
    \left |\Sigma \cap (G\cdot x)\right
    |=r_\Sigma\,\,\mathrm{and}\,\,
    \left |\Lambda \cap (G\cdot y)\right |=r_\Lambda.\end{equation}

    \item[iv):] $G$ acts freely on $M'$.
\end{description}

\begin{defn} Let us set $\tilde Y=:\Upsilon ^{-1}(\Sigma')\subseteq
S^1\times G\times \Lambda$. Let $\pi _\Lambda:S^1\times G\times
\Lambda\rightarrow \Lambda$ be the projection onto the third
summand, and let us set $Y=:\pi _\Lambda(\tilde Y)\subseteq
\Lambda'$.
\end{defn}

\begin{lem} \label{lem:transverse-psi}
Let $\Lambda '=\Lambda \cap X'$. Then there exists an open
neighbourhood $V\subseteq \Lambda$ of $\Lambda'$ such that
$\Upsilon$ is immersive on $S^1\times G\times V$.\end{lem}

\textit{Proof.} This follows from the horizontality of $\Lambda$
and of the $G$-action on $X'$, and from Corollary \ref{cor:xi=0}.

\begin{prop}\label{prop:nice-inverse-image}
Suppose that the hypothesis i), ii) and iii) above are satisfied.
Let $\tilde Y_1,\cdots,\tilde Y_r\subseteq S^1\times G\times
\Lambda'$ be the connected components of $\tilde Y$, and let
$Y_j=:\pi _\Lambda (\tilde Y_j)\subseteq \Lambda'$. Then:
\begin{description}
  \item[i):] for every $j=1,\ldots,r$, there exists
  $h_{j}\in S^1$ such that $$\tilde Y_j\subseteq \{h_{j}\}\times G\times \Lambda';$$
  \item[ii):] every $Y_j$ is
  a submanifold, and the induced map $\pi _j:\tilde Y_j\rightarrow Y_j$ is an unramified
  covering;
  \item[iii):] the $Y_j$'s, with possible repetitions, are the
  connected components of $Y$.
  \end{description}
\end{prop}

\textit{Proof.} i): Suppose that $(h,g,x)\in \tilde Y_j$ for some
$j$, and consider $(a,v,w)\in T_{(h,g,x)}\widetilde{Y}_j$. Since
$\Upsilon (Y_j)\subseteq \Sigma$ and $\Sigma$ is Legendrian, we
conclude that \begin{eqnarray*}0&=&\alpha _{\Upsilon(h,g,x)} \left
(d_{(h,g,x)}\Upsilon(a,v,w)\right )=\alpha _{\Upsilon(h,g,x)}\left
( a\,\frac{\partial}{\partial
\theta}+d_{(h,g,x)}\Upsilon(0,v,w)\right
)\\
&=&a+\alpha _{\Upsilon(h,g,x)}\left (
d_{(h,g,x)}\Upsilon(0,v,w)\right )=a.
\end{eqnarray*} The latter equality follows from the horizontality of
$\Lambda$ and of the $G$-action on $X'$. Since $\tilde Y_j$ is
connected, the statement follows.

ii) and iii): Let $\widetilde{\pi} _j :\tilde Y_j\rightarrow
\Lambda'$ be the projection. If $\widetilde{\pi} _j$ is not an
immersion, by part i) there exists $(h_{j},g,x)\in \tilde Y_j$ and
a tangent vector of the form $(0,v,0)\in T_{(h_{j},g,x)}\tilde
Y_j$, for some $0\neq v\in T_gG$. By Lemma
\ref{lem:transverse-psi},
$$0\neq d_{(h_{j},g,x)}\Upsilon\big ((0,v,0)\big )\in
\left [T_{\Upsilon(h_{j},g, x)}\Sigma\right ]\cap \left
[T_{\Upsilon(h_{j},g, x)}\left (G\cdot \Upsilon(h_{j},g, x)\right
)\right ],$$ against Corollary \ref{cor:xi=0}.

Suppose now that $h\in \{h_1,\ldots,h_r\}\subseteq S^1$. Let
$$Y^{(h)}=:\bigcup _{h_j=h}\tilde Y_j.$$ Suppose that $y\in \pi
_\Lambda ( Y^{(h)})$; there are as many inverse images of $y$ in
$Y^{(h)}$ as there are group elements $g\in G$ such that
$(h,g)\cdot y=g\cdot (h\cdot y)\in \Sigma$; in other words,
\begin{equation}\label{eqn:constant-inv-image}
\left |Y^{(h)}\cap \pi _\Lambda ^{-1}(y)\right |=\left | \Sigma
\cap \big (G\cdot h\cdot y)\right |=d_\Sigma. \end{equation} On
the other hand, in an immersion with compact domain the number of
points in a inverse image can only jump up. Therefore, given
(\ref{eqn:constant-inv-image}) the cardinality of a fibre has to
constant for each map $\tilde Y_j\rightarrow Y_j$, for every
connected component $\tilde Y_{j}$ of $Y^{(h)}$. If on the other
$f:X\rightarrow Y$ is an immersion, and $\left |f^{-1}(y)\right |$
is constant for every $y\in f(X)$, then $f(X)$ is a manifold and
the induced map $X\rightarrow f(X)$ is an unramified covering.

\bigskip

We note in passing that by the same argument, and given the
symmetry of our hypothesis on $\Lambda$ and $\Sigma$, we also
have:

\begin{prop}\label{prop:nice-inverse-image-sigma}
For every $j=1,\ldots,r$, let $\Sigma _j=\Upsilon (\tilde Y_j)$.
Then the $\Sigma _j$'s are disjoint manifolds, and are the
connected components of $\Upsilon (\tilde Y)$. The induced map
$\tilde Y\rightarrow \Upsilon (\tilde Y)$ is an unramified
covering.\end{prop}

\begin{defn} \label{defn:pi_iandsigma}
For every $j=1,\ldots,r$, set $\mathrm{c}_j=:n-\dim (Y_j)$ and let
$\mathrm{d}_j$ be the degree of the unramified cover $\pi
_j:\tilde Y_j\rightarrow Y_j$. Thus, for every $y\in Y_j$ there
exist distinct $s_{1j}(y),\ldots,s_{\mathrm{d}_j,j}(y)\in G$ such
that $(h_j,s_{ij}(y),y)\in \tilde Y_j$, and therefore
$(h_j,s_{ij}(y))\cdot y\in \Sigma'$, $i=1,\ldots,\mathrm{d}_j$.
Locally on $Y_j$ near $y$ we may think of $s_{ij}$'s as $G$-valued
smooth maps. The $s_{ij}$'s are not globally well-defined as
smooth maps $Y_j\rightarrow G$; nonetheless, collectively they do
define a smooth map from $Y_j$ to the appropriate symmetrized
product of $G$.
\end{defn}

Now let $U_j\subseteq \Lambda$ be some arbitrarily small tubular
neighbourhood of the submanifold $Y_j\subseteq \Lambda$. Since
$v_{k,\varpi}=O(k^{-\infty})$ away from $\Sigma'$, in view of
(\ref{eqn:reproducing}) we have
\begin{equation}\label{eqn:reproducing-asympt}
(u_{k,\varpi},v_{k,\varpi})_{L^2(X)}\sim \sum _{j=1}^r\int
_{U_j}\, f_\lambda \cdot
\overline{v_{k,\varpi}}\,\mathrm{dens}_\Lambda.\end{equation}

\begin{rem}\label{rem:not-all-distinct}
Since the $Y_j$'s are not necessarily all distinct,
(\ref{eqn:reproducing-asympt}) is not literally true. However, to
avoid making our exposition too heavy, we shall be slightly vague
on this; we shall thus act as the $Y_j$ were all disjoint. In the
following computations, each summand in
(\ref{eqn:reproducing-asympt}) will split as the sum of various
other contributions, and we shall not sum the same contribution
twice.
\end{rem}

Suppose $1\le j\le r$. For every $y\in Y_j$, we may find an open
neighbourhood $y\in S\subseteq Y_j$ which is uniformly covered by
$\tilde \pi _j$, meaning that $\tilde \pi _j^{-1}(S)=\bigcup
_{i=1}^{d_j}\tilde S_{i}\subseteq \tilde Y_j$, a disjoint union
where each $\tilde S_{i}$ projects diffeomorphically onto $S$
under $\tilde \pi _j$, and $(h_j,s_{ij}(y),y)\in \tilde S_i$.

Perhaps after restricting $S$, by Lemma \ref{lem:transverse-psi}
we may further assume that for each $i$ the map induced by
$\Upsilon$, $(h_j,s_{ij}(y),y)\mapsto (h_j,s_{ij}(y))\cdot y$ is a
diffeomorphism onto its image,

%%%%%%%%%%%%%%%%%%
%%%%%%%%%%%%%%%%%%%%
%%%%%%%%%%%%%%%%%%
\begin{equation}\label{eqn:S_itildeS_i}
\tilde S_i\cong \widehat{S}_i=:\Upsilon (\tilde S_i)\subseteq
\Sigma.
\end{equation}
%%%%%%%%%%%%%%%%%%%%%%%
%%%%%%%%%%%%%%%%%%%%%
%%%%%%%%%%%%%%%%%%%%

We may then find a finite open cover $\{S_{ja}\}_{a\in
\mathcal{A}}$ of $Y_j$ with the following properties:

\noindent i): each $S_{ja}$ is the domain of a coordinate chart,
say
$$
R_{ja}=(r_1,\ldots,r_{\mathrm{n}-\mathrm{c}_j}):S_{ja}\rightarrow
B_{\mathrm{n}-\mathrm{c}_j}(0,\epsilon)\subseteq
\mathbb{R}^{\mathrm{n}-\mathrm{c}_j},$$ for some $\epsilon
>0$;

\noindent ii): each $S_{ja}$ is uniformly covered by $\tilde \pi
_j$, and $\tilde \pi _j^{-1}(S_{ja})=\bigcup _{i=1}^{d_j}\tilde
S_{ija}$ is a disjoint union, where for each $i,\,a$
\begin{equation}\label{eqn:Sijatilde}
\tilde S_{ija}=:\{(h_j,s_{ij}(y),y):y\in S_{ja}\}\subseteq \tilde
Y_j;\end{equation}

\noindent iii): $\Upsilon$ induces a diffeomorphism $\tilde
S_{ija}\cong \widehat{S}_{ija}=:\Upsilon (\tilde S_{ija})\subseteq
\Sigma $ for every $i,a$;

\noindent iv): for every $i,a$ there exist an open neighbourhood
$T_{ija}\subseteq X$ of $\widehat{S}_{ija}$, and a smooth map
$\kappa =\kappa_{ija}:\widehat{S}_{ija}\times T_{ija}\rightarrow
B_{2\mathrm{n}+1}(0,\epsilon)$, such that for every $y\in
\widehat{S}_{ija}$ the partial function $\kappa
_{y}:T_{ija}\rightarrow B_{2\mathrm{n}+1}(0,\epsilon)$ is a
Heisenberg chart adapted to $\Sigma$ at $y$ (Lemma
\ref{lem:variational-version}).

Now recall that for every $j$ we have fixed a tubular
neighbourhood $U_j\subseteq \Lambda$ of $Y_j$; let
$p_j:U_j\rightarrow Y_j$ be the projection, and set
$U_{ja}=:p_j^{-1}(S_{ja})\subseteq U_j$. Thus, $\{U_{ja}\}_{a\in
\mathcal{A} }$ is a finite open cover of $U_j$. By introducing a
partition of unity $\sum _a\varphi _{ja}=1$, we may decompose the
$j$-th summand in (\ref{eqn:reproducing-asympt}) as
\begin{equation}\label{eqn:reproducing-asympt-partitioned}
\int _{U_j}\, f_\lambda \cdot
\overline{v_{k,\varpi}}\,\mathrm{dens}_\Lambda =\sum _a \int
_{U_{ja}}\, \varphi _{ja}\,f_\lambda \cdot
\overline{v_{k,\varpi}}\,\mathrm{dens}_\Lambda .\end{equation} We
are thus reduced to considering the asymptotics of each summand in
(\ref{eqn:reproducing-asympt-partitioned}). Given
(\ref{eqn:Sijatilde}) we may the apply a relative version of the
argument in \S \ref{sctn:transversecase}; rescaling will now be in
the coordinates in $U_{ja}$ which are transversal to $Y_j$.

We now leave it to the reader to verify that, using the local
coordinates $R_{ja}=(r_1,\ldots,r_{\mathrm{n-c}_j})$ on $S_{ja}$,
one obtains an asymptotic expansion
%%%%%%%%%%%%%
%%%%%%%%%%%%%
\begin{eqnarray}\label{eqn:jathtermexpanded}
%%%%%%%%%%%%%
%%%%%%%%%%%%%
\int _{U_{ja}}\, \varphi _{ja}\,f_\lambda \cdot
\overline{v_{k,\varpi}}\,\mathrm{dens}_\Lambda
&\sim&k^{(\mathrm{n-g}-\mathrm{c}_j)/2} \int
_{\mathbb{R}^{\mathrm{n-c}_j}} \rho
_0^{(ja)}(r)\,dr\nonumber\\
&&+\sum _{f\ge
  1}k^{(\mathrm{n-g}-\mathrm{c}_j-f)/2}\,\rho _f^{(ja)}
\end{eqnarray}
where, in view of the asymptotic expansion of Theorem
\ref{thm:main2},
\begin{eqnarray*}
\rho _0^{(ja)}(r)&=&\frac{\dim(V_\varpi)}{V_{\mathrm{eff}}[\pi
(y(r))]}\,
  \frac{1}{\pi ^\mathrm{n}}\,\sqrt{\frac{(2\pi)^{\mathrm{n}+\mathrm{g}}}{2^{\mathrm{g}}}}\,
  h_j^k\,\sum _l\overline{\chi _\varpi \big [s_{lj}(y(r))\big]}\\
  &&\cdot
  \varphi _{ja}\big (y(r)\big
  )\,
  f_\lambda\big (y(r)\big )\,\overline{f_\sigma
 \big [h_j\cdot s_{lj}(y(r))\big]}\,
 \int _{\mathbb{R}^{\mathrm{c}_j}}e^{-S_r(z)-iP_r(z)}\,dz,
\end{eqnarray*}
for quadratic forms $S_r$, $P_r$ on $ \mathbb{R}^{\mathrm{c}_j}$,
with $S_r$ positive definite. In the action-free case, this
becomes:
\begin{eqnarray}\label{eqn:jathtermexpanded}
%%%%%%%%%%%%%
%%%%%%%%%%%%%
\int _{U_{ja}}\, \varphi _{ja}\,f_\lambda \cdot
\overline{v_{k}}\,\mathrm{dens}_\Lambda
&\sim&k^{(\mathrm{n}-\mathrm{c}_j)/2} \int
_{\mathbb{R}^{\mathrm{n-c}_j}} \rho
_0^{(ja)}(r)\,dr\nonumber\\
&&+\sum _{f\ge
  1}k^{(\mathrm{n}-\mathrm{c}_j-f)/2}\,\rho _f^{(ja)}
\end{eqnarray}
where
\begin{eqnarray*}
\rho _0^{(ja)}(r)&=&\frac{\dim(V_\varpi)}{V_{\mathrm{eff}}[\pi
(y(r))]}\,
  \frac{(2\pi)^{\frac{\mathrm{n}}{2}}}{\pi ^\mathrm{n}}\,
  h_j^k\\
  &&\cdot\varphi _{ja}\big (y(r)\big
  )\,
  f_\lambda\big (y(r)\big )\,\overline{f_\sigma
 \big (h_j\cdot y(r)\big)}\,T_{ja} \big(y(r)\big),
\end{eqnarray*}
\begin{eqnarray*}
T_{ja} \big(y(r)\big)&=:& \imath_{J_{h_j\cdot y(r)}}\left
(T_{h_j\cdot y(r)}\cdot \Lambda _j,T_{h_j\cdot y(r)}\Sigma \right
)^{-1}\cdot \int
_{\mathbb{R}^{\mathrm{c}_j}}e^{-\|p\|^2+ip^tZ_{h_j\cdot
y(r')}p}\,dr',
\end{eqnarray*}
$Z_r$ being an appropriate $ \mathrm{c}_j\times \mathrm{c}_j$
symmetric matrix.


\begin{thebibliography}{Dillo99}

\bibitem[BW]{bw} S. Bates, A. Weinstein, {\em
Lectures on the geometry of quantization}, Berkeley Mathematics
Lecture Notes \textbf{8}, AMS 1997

\bibitem[BSZ]{bsz} P. Bleher, B. Shiffman, S. Zelditch, {\em
Universality and scaling of correlations between zeros on complex
manifolds}, Invent. Math. \textbf{142} (2000), 351--395

\bibitem[BPU]{bpu}
D. Borthwick, T. Paul, A. Uribe, {\em Legendrian distributions
with applications to relative Poincar\'{e} series}, Invent. Math.
\textbf{122} (1995), no. 2, 359--402

\bibitem[BS]{bs} L. Boutet de Monvel, J. Sj\"ostrand,
{\em Sur la singularit\'e des noyaux de Bergman et de Szeg\"o},
Ast\'erisque \textbf{34-35} (1976), 123--164

\bibitem[BG]{burnsg} D. Burns, V. Guillemin, {\em Potential
functions and actions of tori on K\"{a}hler manifolds}, Comm.
Anal. Geom. \textbf{12} (2004), no. 1-2, 281--303

\bibitem[Di]{dixmier} J. Dixmier, {\em Les $C^*$-algebras et leurs
r\'{e}pr\'{e}sentations}, Gauthier-Villars Paris (1964)

\bibitem[Du]{duist} J. J. Duistermaat,
{\em Fourier integral operators}, Birkh\"{a}user Boston 1996

\bibitem[Ge]{gei} H. Geiges,
{\em Contact Geometry}, in Handbook of Differential Geometry
\textbf{2}, F.J.E. Dillen and L.C.A. Verstraelen, \textit{eds},
North Holland, Amsterdam (2006), 325--382

\bibitem[GT]{gt} A. L. Gorodentsev, A. N. Tyurin, {\em Abelian Lagrangian algebraic
geometry}, Izv. Ross. Akad. Nauk Ser. Mat. \textbf{65:3} (2001),
15--50; English transl., Izv. Math. \textbf{65} (2001), 437--467

\bibitem[GGK]{ggk} V. Guillemin, V. Ginzburg, Y. Karshon
{\em Moment maps, cobordism, and Hamiltonian group actions},
Mathematical Surveys and Monographs \textbf{98}, A.M.S. (2002)

\bibitem[GP]{gp} V. Guillemin, A. Pollack,
{\em Differential topology}, Prentice-Hall, Inc., Englewood
Cliffs, N.J., 1974


\bibitem[GS1]{gs-gq} V. Guillemin, S. Sternberg,
{\em Geometric quantization and multiplicities of group
representations}, Inv. Math. {\bf 67} (1982), 515--538

\bibitem[GS2]{gs-hq} V. Guillemin, S. Sternberg,
{\em Homogeneous quantization and multiplicities of group
representations}, J. Func. Anal. {\bf 47} (1982), 344--380

\bibitem[GS3]{gs-gc} V. Guillemin, S. Sternberg,
{\em The Gelfand-Cetlin system and quantization of the complex
flag manifold}, J. Func. Anal. {\bf 52} (1983), 106--128

\bibitem[H]{hor} L. H\"{o}rmander,
{\em The analysis of partial differential operators I},
Springer-Verlag 1990

\bibitem[K]{kost} B. Kostant,
{\em Quantization and unitary representations. I.
Prequantization}, Lectures in modern analysis and applications,
III (1965), pp. 87--208. Lecture Notes in Math., Vol.
\textbf{170}, Springer, Berlin, 1970

\bibitem[P1]{pao-mm} R. Paoletti, {\em Moment maps and equivariant Szeg\"{o}
kernels}, J. Symplectic Geom. \textbf{2} (2003), no. 1, 133--175

\bibitem[P2]{pao-sq} R. Paoletti, {\em The Szeg\"{o} kernel of a symplectic quotient},
Adv. Math. \textbf{197} (2005), 523--553

\bibitem[STZ]{stz} B. Shiffman, T. Tate, S. Zelditch, {\em Distribution
laws for integrable eigenfunctions}, Ann. Inst. Fourier (Grenoble)
\textbf{54} (2004), no. 5, 1497--1546

\bibitem[SZ]{sz} B. Shiffman, S. Zelditch, {\em Asymptotics of almost
holomorphic sections of ample line bundles on symplectic
manifolds}, J. Reine Angew. Math. {\bf 544} (2002), 181--222

\bibitem[W]{wein} A. Weinstein, {\em Connections of Berry and Hannay type
for moving Lagrangian submanifolds}, Adv. Math. \textbf{82}
(1990), 133--159

\bibitem[Z]{z} S. Zelditch, {\em Szeg\"o kernels and a theorem of Tian},
Int. Math. Res. Not. {\bf 6} (1998), 317--331





\end{thebibliography}
\end{document}